\documentclass[10pt, a4paper]{amsart}
\usepackage[foot]{amsaddr}

\usepackage{amssymb}
\usepackage{amsthm}

\newtheorem{remark}{Remark}
\usepackage{multirow}
\newcommand{\R}{\mathbb{R}}
\usepackage{amsmath}
\usepackage{amsfonts}
\DeclareMathOperator{\divergence}{div}

\usepackage{tikz}
\usetikzlibrary{calc}

\newcommand{\GammaD}{\Gamma_{\rm D}}
\newcommand{\GammaN}{\Gamma_{\rm N}}
\newcommand{\GammaBC}{\Gamma_{\rm BC}}

\newcommand{\stress}{\boldsymbol \sigma}
\newcommand{\strain}{\boldsymbol \varepsilon}

\makeatletter
\renewcommand{\email}[2][]{%
  \ifx\emails\@empty\relax\else{\g@addto@macro\emails{,\space}}\fi%
  \@ifnotempty{#1}{\g@addto@macro\emails{\mbox{\textrm{(#1)}}\space}}%
  \g@addto@macro\emails{#2}%
}
\makeatother
\begin{document}

\title[A hybrid isogeometric approach on multi-patches]{A hybrid isogeometric approach on multi-patches with applications to   Kirchhoff plates and  eigenvalue problems}
\thanks{
TH, BW, and LW would like to gratefully acknowledge the funds provided by the Deutsche Forschungsgemeinschaft under the
contract/grant numbers: WO 671/11-1, WO 671/13-2 and  WO 671/15-1 (within the Priority Programme SPP 1748).
AR would like to gratefully acknowledge the funds provided by Fondazione Cariplo - Regione Lombardia through the project ``Verso nuovi strumenti di simulazione super veloci ed accurati basati sull'analisi isogeometrica'', within the program RST - rafforzamento. The authors would like to thank Giancarlo Sangalli (University of Pavia) for some fruitful discussions on the topic of this work.
}

 \author{Thomas~Horger}
 \author{Alessandro~Reali}
 \author{Barbara~Wohlmuth}
 \author{Linus~Wunderlich}

 \email[T.~Horger]{horger@ma.tum.de}
 \email[A.~Reali]{alessandro.reali@unipv.it}
 \email[B.~Wohlmuth]{wohlmuth@ma.tum.de}
 \email[L.~Wunderlich, Corresponding author]{linus.wunderlich@ma.tum.de}

 \address[T.~Horger, B.~Wohlmuth, L.~Wunderlich]{Institute for Numerical Mathematics, Technische Universit\"at M\"unchen, Boltzmannstra\ss{}e~3,
85748 Garching b. M\"unchen, Germany}

 \address[A.~Reali]{Department of Civil Engineering and Architecture, University of Pavia, via Ferrata~3, 27100 Pavia, Italy}
 \address[A.~Reali]{Institute for Advanced Study, Technische Universit\"at M\"unchen, Lichtenbergstra\ss{}e~2a,
85748 Garching b. M\"unchen, Germany}
\date{}
\begin{abstract}
We present a systematic study on higher-order penalty techniques for isogeometric mortar methods.  In addition to the weak-continuity enforced by a mortar method,  normal derivatives across the interface are penalized. 
The considered applications are  fourth order problems as well as   eigenvalue problems for second and fourth order equations.
The hybrid coupling enables the discretization of fourth order problems in a multi-patch setting as well as a convenient  implementation of natural boundary conditions. 
For second order eigenvalue problems, the pollution of the discrete spectrum - typically referred to as ``outliers'' - can be avoided.

Numerical results illustrate the good behaviour of the proposed method in simple systematic studies as well as more complex multi-patch mapped geometries for linear elasticity and Kirchhoff plates.
\end{abstract}

\maketitle

\section{Introduction}
Isogeometric analysis (IGA)~\cite{hughes:05} is a family of  methods using  highly regular basis functions typical of CAD systems, like non-uniform rational B-splines (NURBS), to construct numerical approximations of partial differential equations (PDEs). The idea of using spline functions for the approximation of PDEs can   be found in earlier works, see, e.g.,~\cite{Hoellig03}, and were extended by the  isoparametric paradigm, with the goal of   simplifying the mesh generation and refinement processes, possibly bridging the gap between CAD and analysis, see also~\cite{beirao:14,nguyen:15}.

In general, when dealing with non-trivial engineering applications, the computational domain is represented by several spline patches and thus efficient techniques to couple  different patches are required. To retain the flexibility of the meshes at the interfaces, mortar methods are a very attractive option, originally introduced for the coupling of non-matching meshes in spectral and finite element methods~\cite{ben_belgacem:99,bernardi:94,Woh01_a}.
While mortar finite element formulations are quite often motivated by the flexibility of domain decomposition techniques or
by the robustness with respect to non-matching meshes in dynamic applications,
IGA leads in a natural way to a multi-patch situation in case of complex geometries,  see, e.g.,~\cite{hesch:12,bletzinger:14,dornisch:14,nguyen:14,popp:18}. A mathematical stability and a priori analysis enlightening the use of different dual spaces can be found  in~\cite{brivadis:15}. 
In this context also higher-order couplings recently gained attention (see, e.g.,  \cite{kiendl:09, kiendl:10} for Kirchoff-Love shells). 
A  discussion of strong $C^1$ couplings in multi-patch settings is given in \cite{bouclier:16}, whereas  weak continuity of the normal stress is realized in
\cite{collin:16}. Alternative higher-order coupling methods based on least-squares techniques were proposed in~\cite{coox:16}.

Here, we investigate the influence of hybrid couplings on   fourth order problems as well as   on the eigenvalue approximation of second order problems. 
In addition to the weak continuity satisfied in terms of a Lagrange multiplier, we apply  a penalty approach for the jump of normal derivatives at the interfaces and Neumann boundaries.  The weights in the penalization terms are selected such that both the condition number growth rate of the algebraic system and the optimal a priori convergence rate are preserved. 

Standard mortar methods enforce weak $H^1$-conformity, which is insufficient to solve fourth order problems, where a stronger coupling, e.g., by using penalty terms, is needed. 
While standard penalty couplings for fourth order problems, where both the jump of the solution and of the normal derivative are penalized, see~\cite{baker:77,mozolevski:03}, pose a sufficient coupling, they are inconvenient to implement. There the consistency terms include third-order derivatives, which are challenging to transform with NURBS geometries.  In contrast, with the hybrid approach only second order terms need to be computed, which are standard and  often already part of  isogeometric software packages, making the hybrid formulation easy to implement. 
In addition, the Neumann penalty presents a flexible way to treat   natural boundary conditions, which in general is a non-trivial task for plate problems. 

Eigenvalue analysis arises in many important applications in science and engineering, where the amount of eigenvalues of interest can be quite different from application to application. For example in vibroacoustics one is typically only interested in the first part of the spectrum, while for explicit dynamics the approximation quality of the entire spectrum is relevant. 
Compared to FEA, it was observed that IGA possesses superior approximation of eigenvalues, see~\cite{cottrell:06, hughes:08, hughes:09, hughes:14}. Further studies show outliers appearing in the case of reduced continuity~\cite{calo:18} and Neumann boundary conditions~\cite{takacs:16,gallistl:17}.
In our tests a major improvement was shown by a  stronger enforcement of the Neumann boundary condition through a penalty approach, which can  recover spectral results closer to the cases with no outliers. We have to note that the use of penalty introduces high unphysical modes that however, being completely unphysical, can be safely removed, e.g., with a low-rank modification technique (see, e.g.,~\cite{hiemstra2017,reali2016}). 
The better results granted by the proposed method may have a significant impact in those dynamics problems where high modes play an important role and in explicit dynamics, yielding a more favourable CFL condition.

The paper is structured as follows: In Section~\ref{sec:problem}, we briefly review the isogeometric mortar discretization and introduce our hybrid mortar variant. 
For fourth order problems a higher order coupling is necessary to achieve solvability in a nonconforming situation and results are shown in Section~\ref{sec:numerics_plate}.
Numerical results illustrate in Section~\ref{sec:numerics} the influence of the penalization for eigenvalue problems, where the higher part of the spectrum can be significantly improved.
A vibroacoustical example presented in Section~\ref{sec:numerics_violin} presents the application to fourth order eigenvalue problems, where we are interested in the lowest eigenvalues and outliers do not play a significant role. 
Finally, in Section~\ref{sec:conclusion},  conclusions are given. 

\section{Hybrid mortar formulation}
\label{sec:problem}
In this section, we introduce the hybrid mortar method, and later show how it can efficiently be applied to PDEs of second and fourth order.
At first, we briefly recapture the basics of isogeometric mortar methods, and for more details we refer to~\cite{brivadis:15}. For the ease of presentation, we restrict ourselves to the two dimensional case. The generalization to one or three dimensions follows the same lines.

\subsection{Standard mortar coupling}
Let $\Omega \subset \R^2$ be a bounded domain with $\bar{\Gamma}_{\rm D} \cup \bar{\Gamma}_{\rm N}=\partial \Omega$ and $\GammaD \cap \GammaN=\emptyset$. 
Let the domain $\Omega$ be decomposed into $K$ non-overlapping subdomains $\Omega_k$, i.e.,
\[
\overline{\Omega} = \bigcup_{k=1}^K \overline{\Omega}_k, \text{ and } \Omega_i \cap \Omega_j = \emptyset \text{ for } i \neq j.
\]
Here, we limit our presentation to the basic isogeometric concepts and notations used throughout the paper and refer to~\cite{hughes:09,bazilevs:06,piegl:97,Schumaker:07} for more details.
Each of the subdomains is a NURBS geometry, i.e., there exists a NURBS parametrization $\mathbf F_k$ mapping from the parametric space $\widehat \Omega = (0,1)^d$ to $\Omega_k$ based on an open knot vector $\mathbf{\Xi}_k$ and a degree $p$. 
Let us consider  a PDE of order $2n$ and $p\geq n$. 

We set $N^p(\mathbf{\Xi}_k)$ as the multivariate NURBS space (associated to $\Omega_k$) in the parametric domain, with the standard NURBS basis functions $\widehat{N}_{k,\mathbf{i}}^p$.
For a set of control points $\mathbf{C}_{k,\mathbf{i}} \in \mathbb{R}^d$, $\mathbf{i} \in \mathbf{I}$, we define a parametrization of a NURBS surface  as a linear combination of the basis functions and  control points
\begin{equation*}
	\mathbf{F}_k({\boldsymbol \zeta})=\displaystyle \sum_{\mathbf{i} \in \mathbf{I}} \mathbf{C}_{k,\mathbf{i}} \,\widehat{N}_{k, \mathbf{i}}^p({\boldsymbol \zeta }),
\end{equation*}
and assume the regularity stated in~\cite[Assumption 3.1]{beirao:14}.

For $1\leq k_1 { < }  k_2 \leq K $, we define the interface as the interior of the intersection of the boundaries, i.e., $\overline{\gamma}_{k_1k_2} = \partial {\Omega}_{k_1} \cap \partial {\Omega}_{k_2}$, where ${\gamma}_{k_1k_2}$ is open. Let the non-empty interfaces be enumerated by $\gamma_l$, $l = 1,\,\ldots,\, L$. 

For each $\Omega_k$, we introduce  $H^n_*(\Omega_k)=\{ v_k \in H^n(\Omega_k),  v_{k |_{\GammaD \cap \partial \Omega_k}}=0\}$, where we use standard Sobolev spaces, as defined in~\cite{grisvard:11}, endowed with their usual norms.  In order to set a global functional framework on $\Omega$, we consider the broken Sobolev space $V= \prod_{k=1}^K H^n_*(\Omega_k)$, endowed with the broken norm $ \| v \|_{V}^2 = \sum_{k=1}^K \| v \|_{H^n(\Omega_k)}^2$.
For any interface $\gamma_l \subset \partial \Omega_{k}$, we define  
{ $H^{-1/2}(\gamma_l)$ to be the dual space of $H^{1/2}_{00}(\gamma_l)$, which is the space of all functions that can be trivially extended (i.e. by zero) on $\partial \Omega_{k} \setminus \gamma_l$  to an element of $H^{1/2}(\partial \Omega_{k})$.}

In the following, we set our non-conforming approximation framework. On each subdomain $\Omega_k$, based on the NURBS parametrization, we introduce the approximation space $V_{k}=\{v_k=\widehat{v}_k \circ \mathbf{F}_k^{-1}, \widehat{v}_k \in N^{p}(\mathbf{\Xi}_k) \}$. On $\Omega$, we define the discrete product space $V_h = \prod_{k=1}^K V_{k} \subset V$, which forms a non-conforming space with respect to $H^n(\Omega)$.

On the skeleton $\Gamma = \bigcup_{l=1}^L\gamma_l$, we define the discrete Lagrange multiplier product space $M_h$ as $M_h = \prod_{l=1}^L M_{l}$. Based on the interface knot vector of one of the adjacent subdomains, $M_{l}$ is the spline space of degree $p$ defined on the interface $\gamma_l$. 
An appropriate local degree reduction performed at the crosspoints guarantees the inf-sup stability of the mortar coupling, see~\cite{brivadis:15} for more details, while preserving optimal order error decay rates.

The  coupling bilinear form
\[
 b({ \tau, v} ) = \sum_{l=1}^L \int_{\gamma_l} \tau [v]_l ~\mathrm{d}\sigma, \]
 where $[\cdot]_l$ denotes the jump over $\gamma_l$, defines the weakly coupled space
 \[
 X_h = \{v_h\in V_h: b(\tau, v_h) = 0, \tau \in M_h\}.
 \]
We note that for second order PDEs ($n=1$), this space is weakly-conforming, while it is still a weakly non-conforming space for higher order PDEs ($n\geq 2$) as there are no restrictions to the normal derivatives across the interface.

\subsection{Hybrid mortar formulation} 
To improve the global continuity, we penalize the jump in the normal derivatives across the interfaces of a multi-patch geometry and the first normal derivative on a Neumann boundary part. To avoid locking, we take the local $L^2$-projection $\pi_h^0$ onto the piecewise constant functions on the (slave) boundary mesh. Due to the non-conformity of $V_h$, we require $m\geq n-1$ and note that for $n>1$ consistency terms are necessary, which will be introduced later. 

The extra penalty term is given as 
\[c_h(u_h,v_h) =C_{\rm BC}\, c_{\rm BC}(u_h, v_h) +  \sum_{l=1}^L \sum_{m=1}^{p-1} C_l^m\, c_l^m (u_h,v_h)
+ C_{\rm CP} \,c_{\rm CP} (u_h, v_h) 
\] with appropriate penalty constants $C_l^m,  C_{\rm BC}, C_{\rm CP} \geq0$ and  problem-dependent boundary terms $c_{\rm BC}$. In the numerical results, we will not distinguish between the different penalty constants and simply refer to them by $C$.  The smooth interface coupling
\[
c_l^m(u_h,v_h) = \int_{\gamma_l} h_s^{2(m-n)-1} \,  \pi_h^0\left( [ \partial_n^m u_h ]_l \right)\, \pi_h^0 \left( [\partial_n^m v_h]_l \right) \, \rm d \sigma,
\]
and the boundary penalty term (where the penalty boundary part $\GammaBC \subset \partial \Omega$  denotes the part of the boundary, where $\partial_{\bf n} u = 0$)
\[
c_{\rm BC} (u_h, v_h)  = \int_{\Gamma_{\rm BC} }  h^{2(m-n)-1} \, \pi_h^0\left( \partial_n u \right)\, \pi_h^0\left( \partial_n v\right) \, \rm d \sigma,
\]
are properly weighted with the local mesh-size $h_s$ on the slave side. The index~CP in the bilinear form $c_{\rm CP}$ refers to contributions from the crosspoints. At the end points of each interface and each corner of the penalized boundary, we introduce additional point evaluations. More precisely, for each interface $\gamma_l$, we add  the term \[
\sum_{\bar {\bf x} \in \partial \gamma_l} \sum_{m=1}^{p-1}  h_s^{2m-2n} [\partial_n^m u_h]_l(\bar {\bf x}) [\partial_n^m v_h]_l(\bar {\bf x}),
\] 
while an analogous term is added on each corner of the penalty boundary  part~$\Gamma_{\rm BC}$.

As already mentioned, weak $C^1$-continuity can be imposed in terms of an additional Lagrange multiplier approach. However the choice of the discrete Lagrange multiplier space is delicate, since uniform stability has to be guaranteed. Moreover, it   involves a careful handling of the resulting algebraic system. Hence, we herein propose an alternative penalty approach.
\section{Application to Kirchhoff plate problems}
\label{sec:numerics_plate}
In this section, we present the   approximation of fourth order problems in multi-patch situations. We solve the bilaplace equation for  clamped Kirchhoff plates
\begin{align*}
\Delta\Delta u &= f\quad \text{ in } \Omega, \\
u&=0 \quad \text { on } \ \partial \Omega, \\
\partial_{\bf n} u &= 0 \quad \text{ on } \ \partial \Omega.
\end{align*}

\subsection{Consistency terms}
\label{sec:consistency_terms}
While for second  order problems, the penalty approach was used only to enforce additional smoothness, for fourth order problems it is necessary to enforce conformity. Hence, it must be used in a Nitsche-type version with additional consistency terms.

We here adapt the symmetric $C^0$ interior penalty Galerkin method for biharmonic formulation of~\cite{brenner:05,brenner:15}: 
Find $ (u_h, \widehat{\tau}_h) \in V_h \times M_h$,  such that 
 \begin{align*}
		a_h^{\rm{bi}}(u_h, v_h)+ b({ \widehat{\tau}_h, v_h }) &= f(v_h), \quad   v_h \in V_h,\\
		b({ \tau_h, u_h }) &= 0, \quad  \tau_h \in M_h,
\end{align*}
with
\begin{align*}
a_h^{\rm{bi}}(u_h, v_h)  &= 
\sum_{k=1}^K \int_{\Omega_k} 
 \operatorname{D}^2 u :  \operatorname{D}^2 v ~ \mathrm{d}{\bf x}
+ c_h(u_h,v_h) \\
&\quad + \sum_{l=1}^L \int_{\gamma_l} \{\partial_{\bf nn} u \}[\partial_{\bf n} v] +[\partial_{\bf n} u] \{\partial_{\bf nn} v \}  \rm d \sigma
\\&\quad + \int_{\GammaD} \partial_{\bf nn} u\, \partial_{\bf n} v  + \partial_{\bf n} u\,\partial_{\bf nn} v ~\rm d \sigma,
\end{align*}
with the Hessian $\operatorname{D}^2v\colon \Omega\rightarrow\mathbb{R}^{2\times 2}$ and $f(v) = \sum_{k=1}^K \int_{\Omega_k}  f v ~\mathrm{d}{\bf x}$.
Since the boundary conditions include the normal derivative, the whole boundary is included in the penalty, i.e., $\GammaBC =  \partial \Omega$. 

We note that the expected convergence rate in the $L^2$ norm is different from the second order case in the lowest order case $p=2$. While Aubin-Nitsche-trick for conforming approximations can be used as for $H^1$-conforming problems, even with optimal dual regularity, the dual approximation order is insufficient to prove the convergence order $h^3$:
\begin{align*}
\left| u-u_h\right|_0^2 &= (u-u_h, u-u_h)_0 = a(w, u-u_h) =  a(w-w_h, u-u_h) \\&\leq \left|   w-w_h \right|_2 \left|   u-u_h \right|_2  \leq c h^2 \|w\|_3 \|u\|_3 \leq c h^2 \left| u-u_h\right|_0 \|u\|_3,
\end{align*}
with the dual solution $w\in H^4(\Omega)$ and its best-approximation $w_h$. 
For $p\geq 3$ the approximation order is sufficient to show optimal order $h^{p+1}$ convergence in the $L^2(\Omega)$ norm. 

\subsection{Multi-patch convergence on a square} \label{subsec:numerics_plate_square}
%At first, let us consider a manufactured solution with homogeneous boundary conditions
As a first numerical test, let us consider a problem with a manufactured solution
 in order to observe the optimality of the method.
Note that all numerical simulations in this article are based on the isogeometric Matlab toolbox GeoPDEs~\cite{geopdes:11,geopdes:16}.

We consider $\Omega = (0,2)\times(0,1)$ and the manufactured solution
\[
u(x,y) = \left(1 - \cos(\pi/2 x) - x + \sin(\pi x)/\pi\right)\left(1-\cos(2\pi y)\right).
\]

\begin{figure}
\includegraphics[height = 11em]{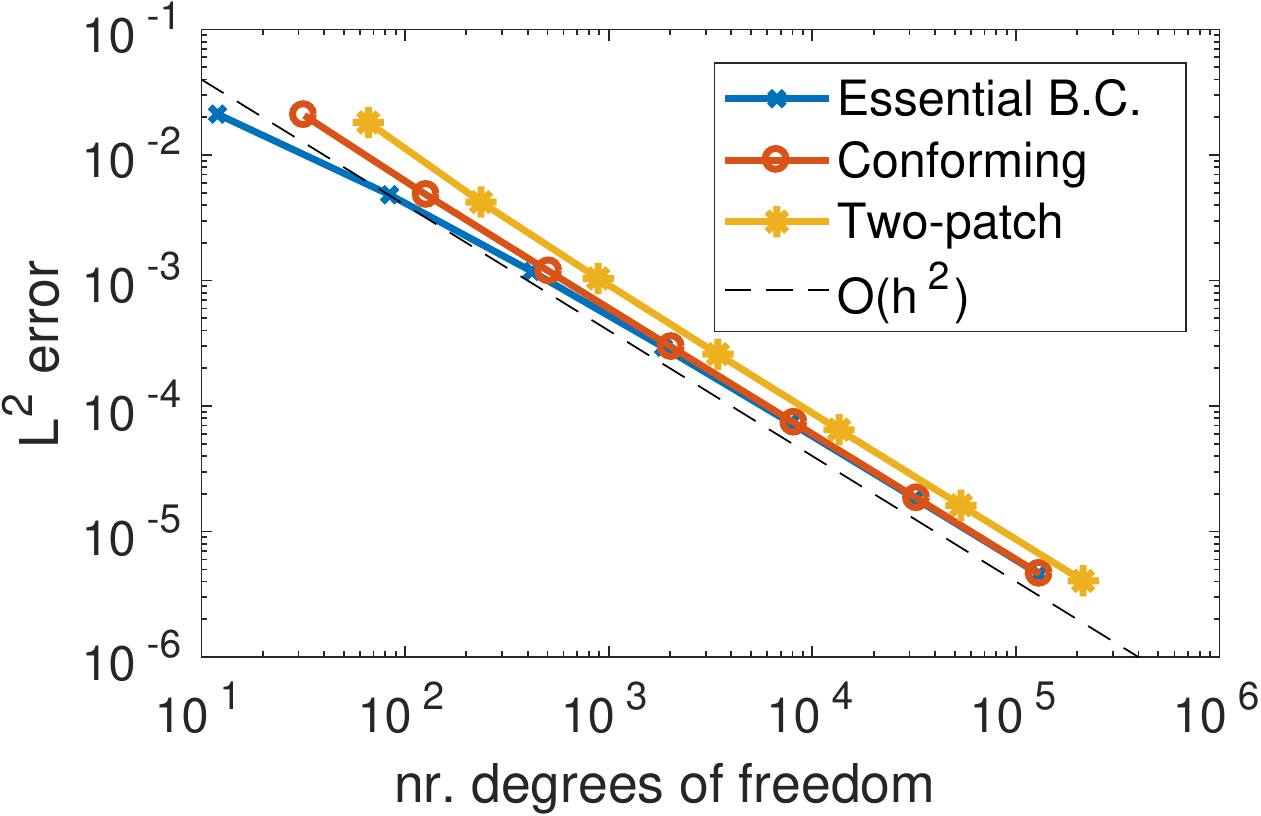}\hfill
\includegraphics[height = 11em]{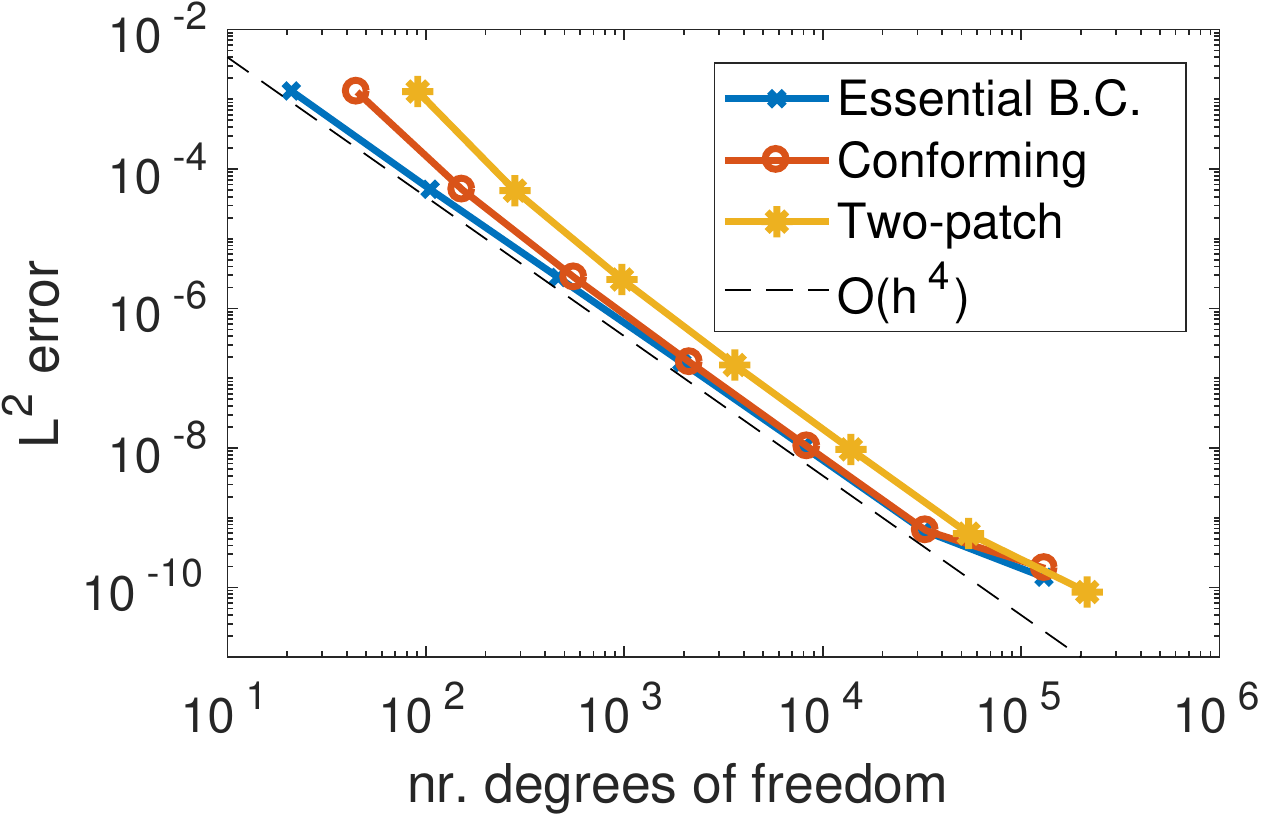}
\caption{$L^2$-convergence for essential boundary conditions, a one-patch penalty and a two-patch hybrid penalty setting with the expected order of convergence. Left: $p=2$. Right: $p=3$.}
\label{fig:plate_squares_convergence}
\end{figure}

We compare three cases: two single-patch settings and a two-patch setting with a non-matching interface. In the first case, both $u=0$ and $\partial_{\bf n} u = 0$ are implemented as essential boundary conditions while in the second case we apply the penalty method on the boundary. The two-patch setting includes penalty terms on the interface as well as the boundary. 
 The convergence in the $L^2$ norm for $C=100$ is shown in Figure~\ref{fig:plate_squares_convergence}.
 We see almost identical optimal error values on the same mesh level and all cases with only the number of degrees of freedom varying. As the boundary values are fixed, the number of degrees of freedom  for the essential boundary conditions are the smallest, but the difference with respect to the one-patch penalty case decreases. In the two-patch setting the number of degrees of freedom is the largest, mainly due to the artificially constructed nonconforming mesh. However, we note that this is a quite artificial setting. In many cases, a single-patch discretization is not possible and a higher-order coupling is essential as seen in the next example. 

In Table~\ref{tab:plate_squares_var_penalty}, the $L^2$ error for the two-patch setting is shown for different values of the penalty constants. We observe robustness within a wide range of penalty values.  

\begin{table}\footnotesize
\hspace*{-2.5em}
\begin{tabular}{r|c|c|c|c|c|c|c}
ndof & 91& 281& 973& 3\,605& 13\,861& 54\,341& 215\,173\\ \hline
\vphantom{ $\Big($ }$C=1$\hphantom{$1$}&  $0.0039$ & $1.53\cdot 10^{-4}$ & $5.47\cdot 10^{-6}$& $2.60\cdot 10^{-7}$&$1.35\cdot10^{-8}$& $7.85\cdot 10^{-10}$& $7.04\cdot10^{-11}$ \\
\vphantom{ $\Big($ }$C= 10^2$& $0.0013$& $4.88\cdot 10^{-5}$& $2.60\cdot 10^{-6}$& $1.55\cdot 10^{-7}$& $9.51\cdot 10^{-9}$& $5.92\cdot 10^{-10}$& $8.57\cdot 10^{-11}$ \\
\vphantom{ $\Big($ }$C= 10^4$& $0.0013$& $4.98\cdot 10^{-5}$& $2.63\cdot 10^{-6}$& $1.56\cdot 10^{-7}$& $9.54\cdot 10^{-9}$& $5.97\cdot 10^{-10}$& $1.16\cdot 10^{-10}$\\
\vphantom{ $\Big($ }$C=  10^6$&$0.0013$& $4.98\cdot 10^{-5}$&$2.63\cdot 10^-6$ &$1.57\cdot 10^{-7}$& $9.95\cdot 10^{-9}$& $5.69\cdot 10^{-9}$&$1.83\cdot 10^{-8}$   
   \end{tabular}
   \caption{$L^2$ error values for the two-patch setting with $p=3$ and a varying value of the penalty parameter.}
   \label{tab:plate_squares_var_penalty}
   \end{table}

\subsection{Beam with holes}
\begin{figure}\begin{center}
\begin{minipage}{.54\textwidth}
\includegraphics[width=\textwidth]{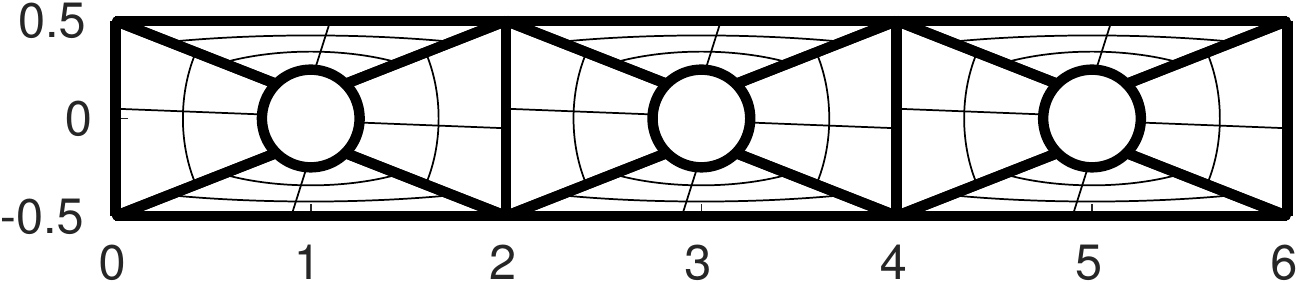}\hfill
\end{minipage}
\begin{minipage}{.45\textwidth}
\includegraphics[width=\textwidth]{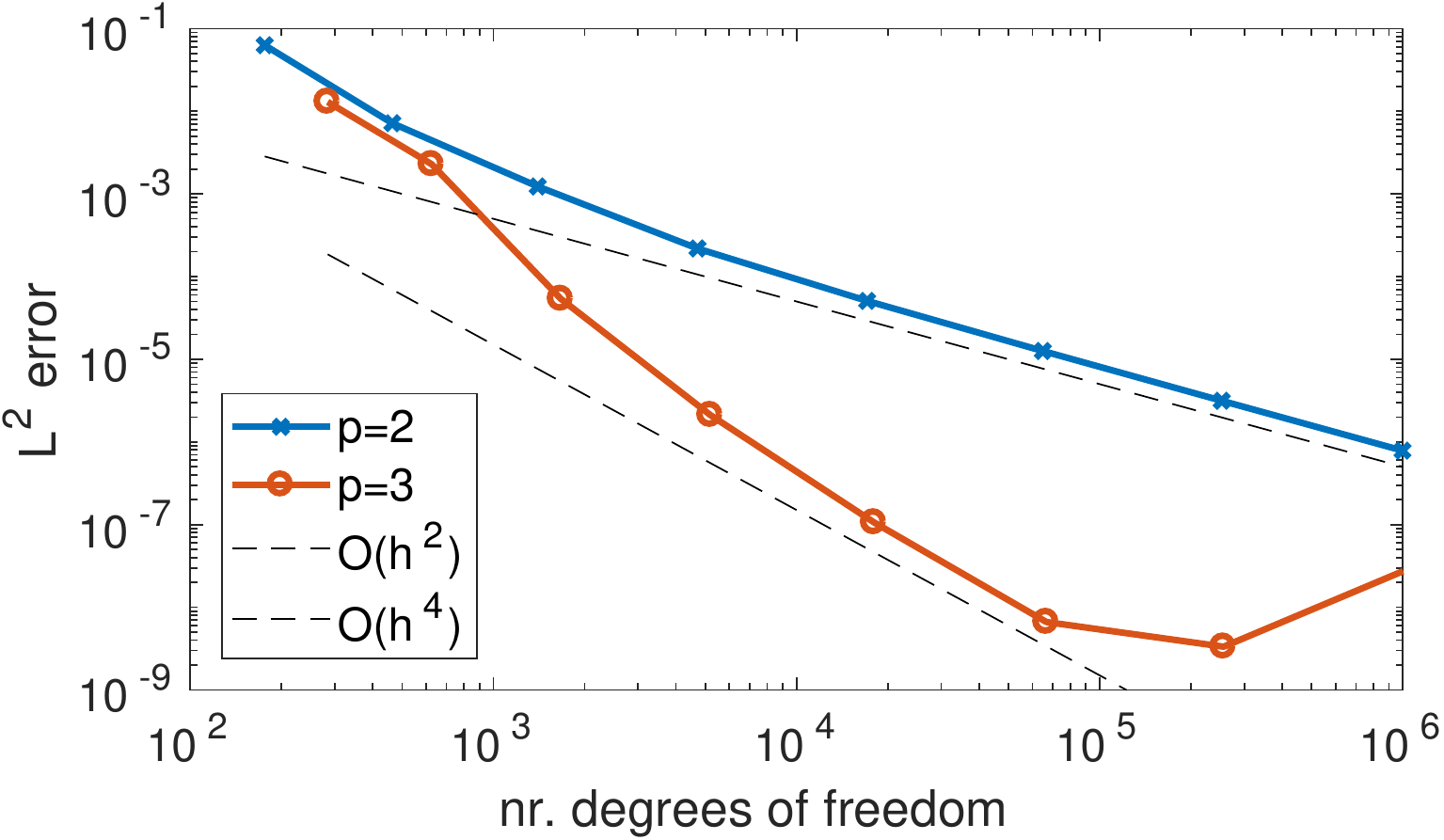}
\end{minipage}
\\[.5em]
\hspace{.5em}
\includegraphics[width=.45\textwidth]{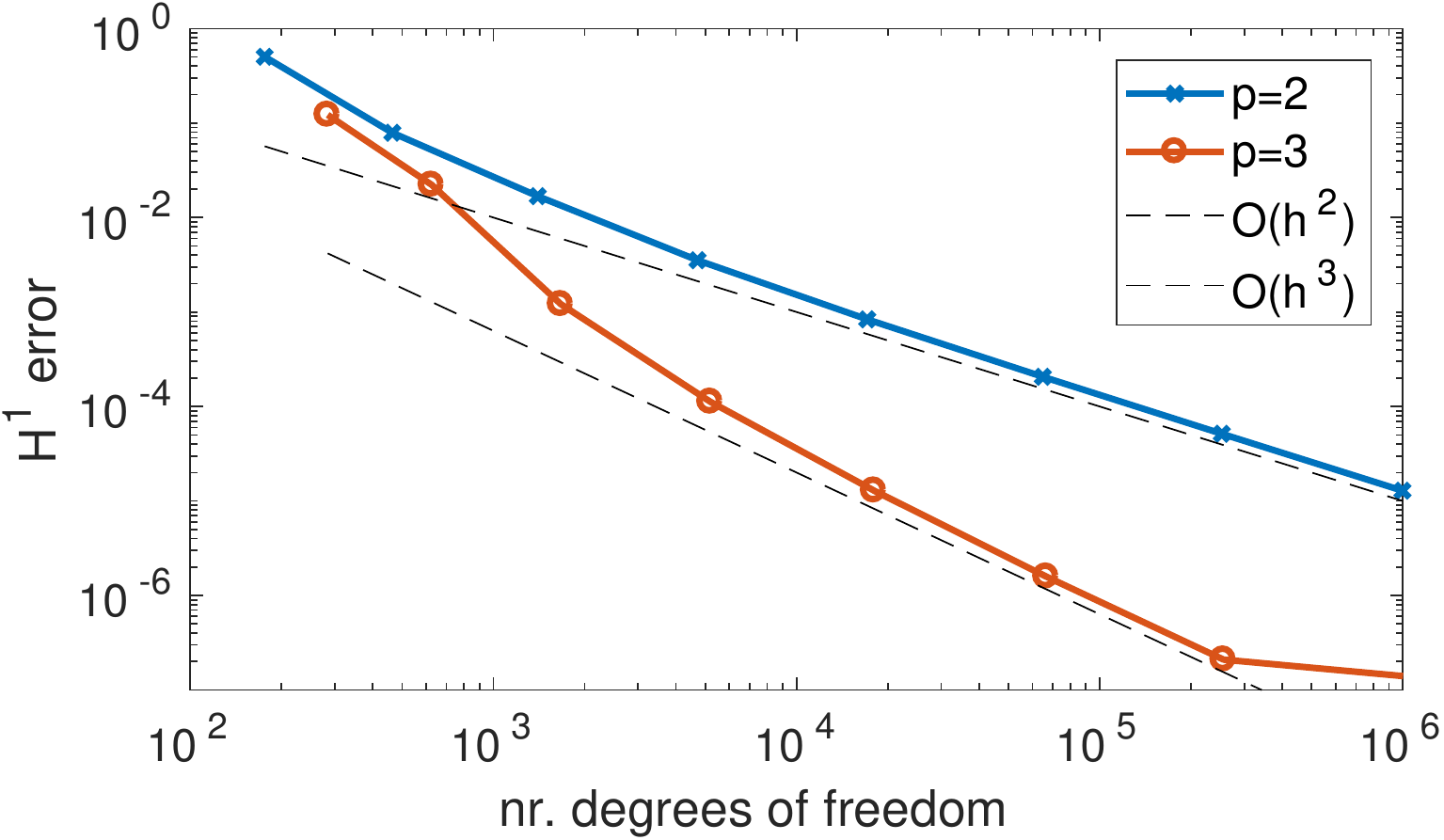}\hfill
\includegraphics[width=.45\textwidth]{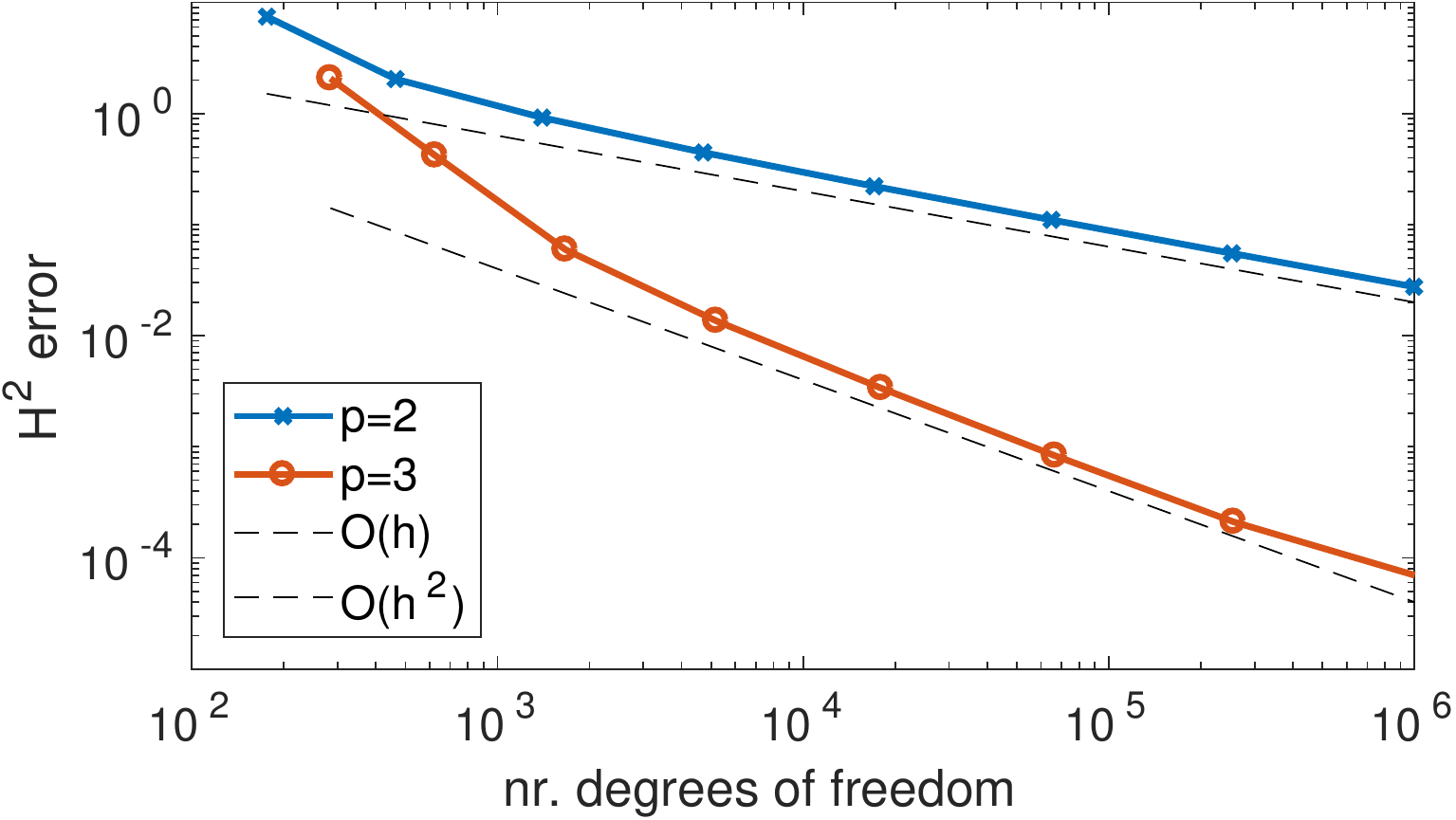}
\end{center}
\caption{Top left: Initial mesh for the beam with circular cut-outs and its decomposition into patches; Convergence for the plate equation on the beam with holes and the expected convergence order for $p=2,3$. Top right: $L^2$ error. Bottom left: $H^1$ error. Bottom right: $H^2$ error.}
\label{fig:beam_with_holes}
\end{figure}

We then consider a beam, with three circular cut-outs, as depicted in Figure~\ref{fig:beam_with_holes}  and observe the convergence for the manufactured solution $u(x,y) = \sin(x) \cos(\pi y)$. 
We note that, although atypical in practice, the manufactured solution has nonhomogeneous boundary conditions for practical reasons. However, this sheds a light on the flexibility in the treatment of the boundary condition, which can be used for practical cases, e.g., when the solution is not restricted, but the normal derivative is.

The resulting errors in the $L^2$, $H^1$ and $H^2$ norms are shown in Figure~\ref{fig:beam_with_holes} for a penalty value of $C=10^3$. We observe  the theoretically expected convergence order, which equals to the  best-approximation order, except for the $L^2$ norm in the quadratic case. 
We also observe, that the error stagnates at a level significantly higher than machine precision, which results from the higher condition number of fourth order equations. 

\section{Application to second order eigenvalue problems}

\label{sec:numerics}
In this section, we present the effect on eigenvalue approximations for second order equations ($n=1$). 
We consider the following Laplace eigenvalue problem with Dirichlet and Neumann boundary  conditions on $\GammaD$ and $\GammaN $, respectively:
\begin{align*} 
		-\Delta u &= \lambda u \quad \text{ in } \Omega,\\
		    u&=0 \quad \text{ on } \GammaD, \\
		    \partial_n u&=0 \quad \text{ on } \GammaN.
\end{align*}
Since the Neumann boundary condition sets the normal derivative, we have $\GammaBC = \GammaN$.
We consider the bilinear forms $a\colon V \times V \rightarrow \mathbb{R}$ and  
$m\colon V \times V \rightarrow \mathbb{R}$, such that  
 \begin{gather*}
 a(u,v) = \sum_{k=1}^K \int_{\Omega_k} 
 \nabla u \cdot  \nabla v ~ \mathrm{d}{\bf x}, \quad
m(u,v) = \sum_{k=1}^K \int_{\Omega_k} u\, v ~ \mathrm{d}{\bf x}
\end{gather*}

 The saddle point formulation  of the isogeometric mortar eigenvalue problem  introduces  $2 \operatorname{dim } M_h $ spurious eigenvalues to the spectrum. We restrict ourselves to the   physical relevant eigenpairs $(\lambda^h, u_h)$. These are characterized by the fact that they are also eigenpairs of the constrained mortar formulation, i.e., they satisfy
\[
a(u_h, v_h) = \lambda^h\, m(u_h, v_h), \quad v_h\in  X_h = \{v_h\in V_h :  b(\tau_h, v_h) = 0, \tau_h \in M_h\}.
\]

In the systematic study of~\cite{gallistl:17}, it was shown that the spectrum of an isogeometric discretization shows  severe outliers in the case of Neumann boundary conditions. The same can be expected for interfaces with $C^0$-regularity, see also~\cite{calo:18}. This motivates us to impose higher-order penalty terms in the formulation:
\begin{align*}
		{a}_h(u_h, v_h)+ b({ \widehat{\tau}_h, v_h }) &= \lambda_h\, m(u_h,v_h), \quad   v_h \in V_h,\\
		b({ \tau_h, u_h }) &= 0, \quad  \tau_h \in M_h,
\end{align*}
with ${a}_h(u_h,v_h) = a(u_h,v_h) + c_h(u_h,v_h)$.

\begin{remark}
For the penalized bilinear form ${a}_h$ broken $H^1$ continuity and for $\left|\Gamma_D\right| > 0$  ellipticity on the kernel of the mortar coupling can   be shown. The ellipticity trivially follows from the ellipticity of $a$, 
while to show continuity, it remains to prove $\left| c_h(u_h, v_h) \right| \leq C \|u_h\|_{V_h} \| v_h \|_{V_h}$. With standard estimates and stability of the $L^2$-projection,  this reduces to an inverse inequality $\| \partial_n^m v_h \|_{L^2(\gamma_l)} \leq C h^{1/2-m} \|v_h\|_{H^1(\Omega_k)}$ for $ m<p$. Then 
 standard trace and inverse inequalities (see~\cite[Theorem 4.2]{bazilevs:06}) yield
 \begin{align*}
 \| \partial_n^m v_h \|_{L^2(\gamma_l)}^2 &\leq C \|  v_h \|_{H^m(\Omega_k )} \|  v_h \|_{H^{m+1}(\Omega_k )}\\
 &\leq   C h^{1-m} \|  v_h \|_{H^1(\Omega_k )} h^{1-(m+1) }\|  v_h \|_{H^1(\Omega_k )}.
 \end{align*}
The  point evaluations can be handled analogously using an inverse inequality between $L^\infty$- and $L^2$-norms.

The a priori analysis of the new hybrid mortar approach can be easily worked out within the abstract framework of non-conforming finite element techniques. The hybrid form allows us to use the Lemma of Strang to show optimal order a priori bounds for right hand side problems, which are required for optimality of the approximation for eigenvalue problems. 
For convenience of the reader, let us sketch the proof for the first order penalty terms and without loss of generality, for a single interface $\gamma$ and no Neumann boundary. 

Let a right hand side $f$ be given and denote the solution to the continuous problem by $u\in H_*^1(\Omega)$. Furthermore, we consider $u_h \in X_h$ the standard mortar discretization and $\widehat{u}_h \in X_h$ the new hybrid solution, which solve
\[
a(u_h,v_h)=f(v_h), \text{ and }
 a_h(\widehat u_h,v_h)=f(v_h),
\]
for each $v_h\in X_h$. 
Since it is well-known~\cite{brivadis:15}, that $u_h$ converges with optimal order, it is sufficient to consider $\|u_h - \widehat u_h \|_{V_h}$ in more detail. Using the coercivity of $a_h$ and a modified Galerkin orthogonality results in 
\begin{align*}
\|u_h - \widehat u_h\|_{V_h}^2 &\leq   a_h(u_h-\widehat u_h, u_h-\widehat u_h) = c_h(u_h,u_h-\widehat u_h)\\
&\leq c \|h^{1/2} [\partial_{\bf n} {u}_h]\|_{L^2(\gamma)} \|h^{1/2} [\partial_{\bf n} (u_h - \widehat u_h)]\|_{L^2(\gamma)}.
\end{align*} 
For the first term, we introduce a suitable best-approximation $w_h\in X_h$ and use $[u] = 0$ for the exact solution:
\begin{align*}
\int_\gamma h [\partial_{\bf n} {u}_h]^2\,\mathrm{d}\sigma & =  
\int_\gamma h [\partial_{\bf n} ({u}_h-w_h)]^2\,\mathrm{d}\sigma+
\int_\gamma h [\partial_{\bf n} ({w}_h-u)]^2\,\mathrm{d}\sigma.
\end{align*}
A polynomial inverse estimate is used for the discrete term, while a local approximation property is used for the second term. Both terms yield an estimate by $ c h^{2s} \|u\|_{H^{s+1}(\Omega)}^2$.

The remaining term $\|h^{1/2} [\partial_{\bf n} (u_h - \widehat u_h)]\|_{L^2(\gamma)}$ can be traced back to the $V_h$ error by polynomial inverse estimates:
\[
\|h^{1/2} [\partial_{\bf n} (u_h - \widehat u_h)]\|_{L^2(\gamma)}\leq
c\| u_h - \widehat u_h\|_{V_h},
\]
which yields the optimal error estimate
$\|u-\widehat u_h\|_{V_h}\leq c h^s \|u\|_{H^{s+1}(\Omega)}$.

We point out that the weights in the penalty term are selected such that the condition number of the algebraic system is still ${\mathcal O} (h^{-2})$. 
\end{remark}

After a systematical one-dimensional investigation, we study a non-trivial multi-patch example in the framework of linear elasticity. 
We compare globally smooth spaces, $C^0$-couplings and the previously introduced higher-order penalty couplings and report the normalized discrete eigenvalue $\lambda^h / \lambda$, which directly relates to the relative error in the eigenvalue since $(\lambda^h - \lambda)/\lambda = \lambda^h/\lambda - 1$. Note that in the case of pure Neumann boundary conditions, the first eigenvalue is zero, so we exclude it from the spectrum.

\subsection{Influence on eigenvectors and the approximation property}

The numerically obtained eigenvalues can be grouped into physical relevant eigenvalues and  unphysical ones induced by the coupling or the boundary. These spurious eigenvalues are infinite for the mortar case and very large in comparison to the physical ones in the penalty case.
In this work we choose to neglect these unphysical modes and only show the physical part of the resulting spectrum. 
To distinguish between physical and unphysical parts of the spectrum, we use a heuristic criterion, namely, $\lambda_{n+1}^h / \lambda_n^h > 100$. 

An ``outlier reduction technique'' based on low-rank modifications, as proposed in \cite{hiemstra2017,reali2016}, conveniently allows to remove unphysical eigenvalues from the spectrum.
Clearly, such a technique does not negatively impact  the approximation properties of the method if the removed modes are actually unphysical and, therefore, do not significantly contribute to the overall response.
In the following, we show with an illustrative example how the high-frequency eigenmodes induced with the $C^0$-coupling contribute to the approximation property of the space and cannot be simply removed. In contrast, with the penalty coupling, the highest modes are indeed unphysical and can be  then safely removed from the space.

\begin{figure}\begin{center}
\includegraphics[width=.44\textwidth]{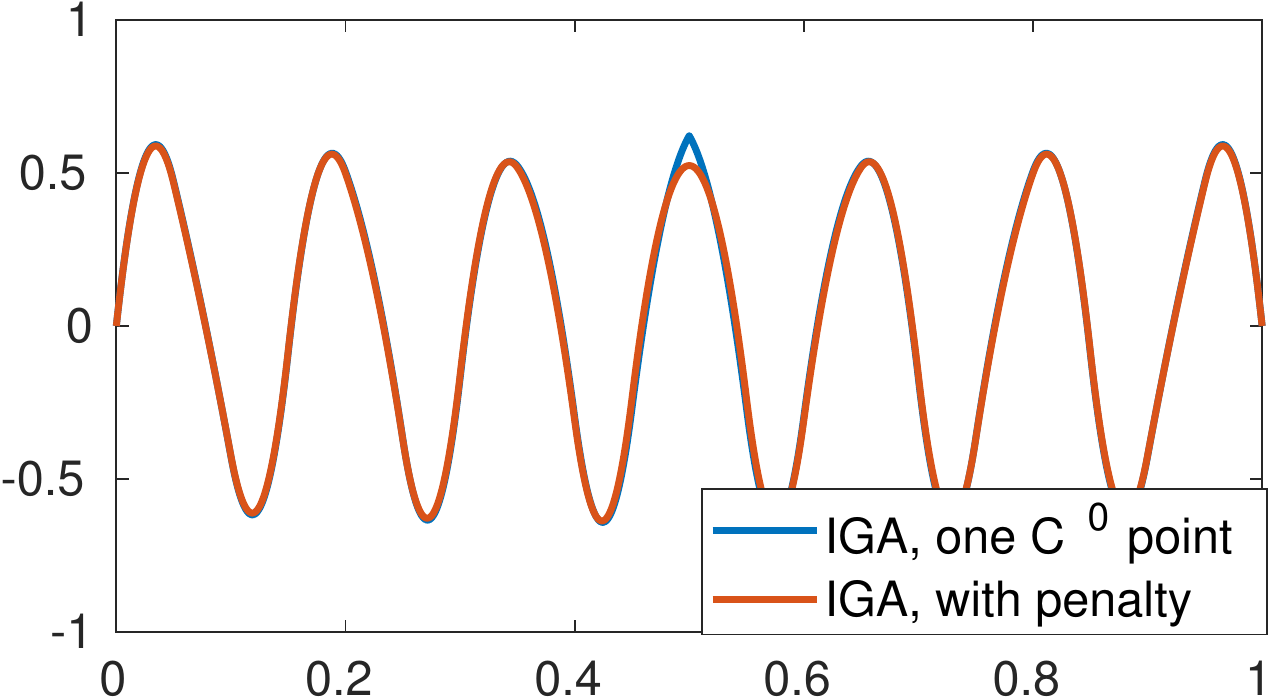}\hspace{1em}
\includegraphics[width=.44\textwidth]{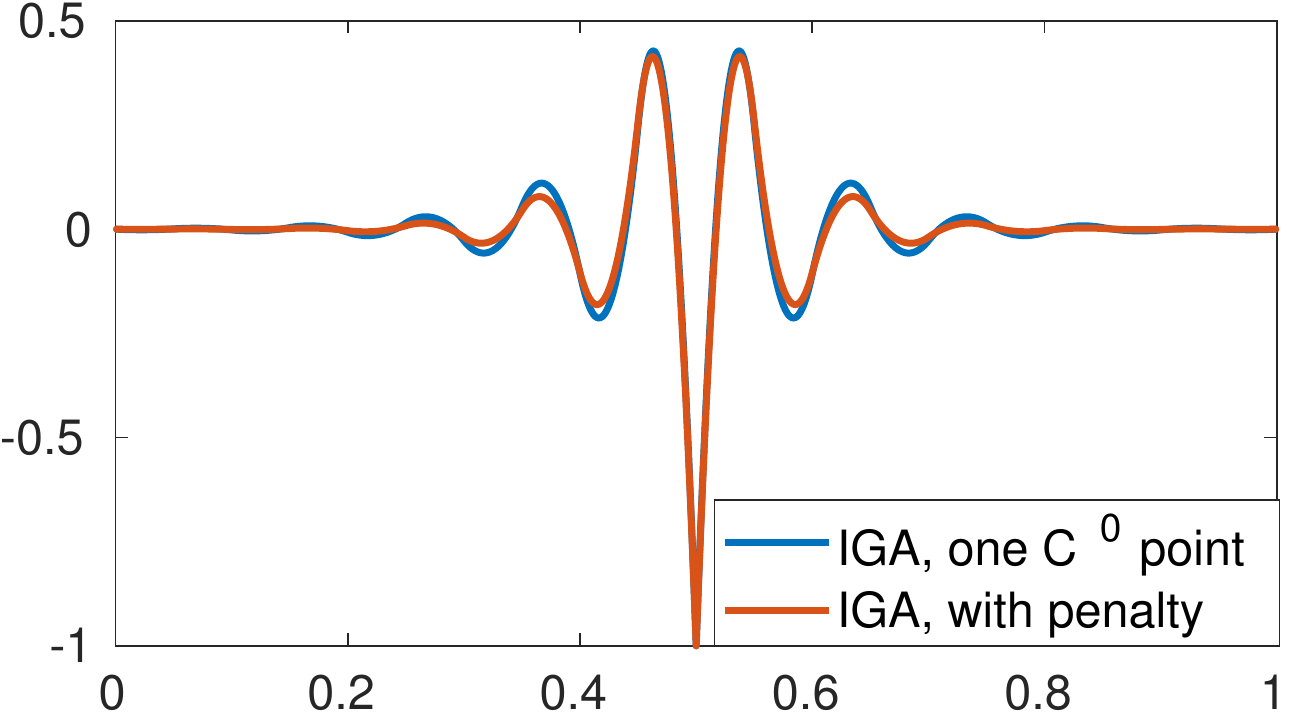}\\
\caption{Discrete eigenvectors for a penalty value of 100. Left: 13th. Right: 21st.}
\label{fig:ex_simple_eigenvectors}
\end{center}\end{figure}

We consider $p=2$ in one dimension on a uniform grid with 21 degrees of freedom once Dirichlet boundary conditions are imposed. Figure~\ref{fig:ex_simple_eigenvectors} compares two selected eigenvectors obtained in the standard case with a $C^0$ point and in the penalty case. While the last eigenvector looks similar in both cases (Figure~\ref{fig:ex_simple_eigenvectors}, right), some of the first 20 ones for the $C^0$ case are different from those for the penalty, as they are non-smooth (see, e.g., Figure~\ref{fig:ex_simple_eigenvectors}, left). 

\begin{figure}\begin{center}
\includegraphics[width=.49\textwidth]{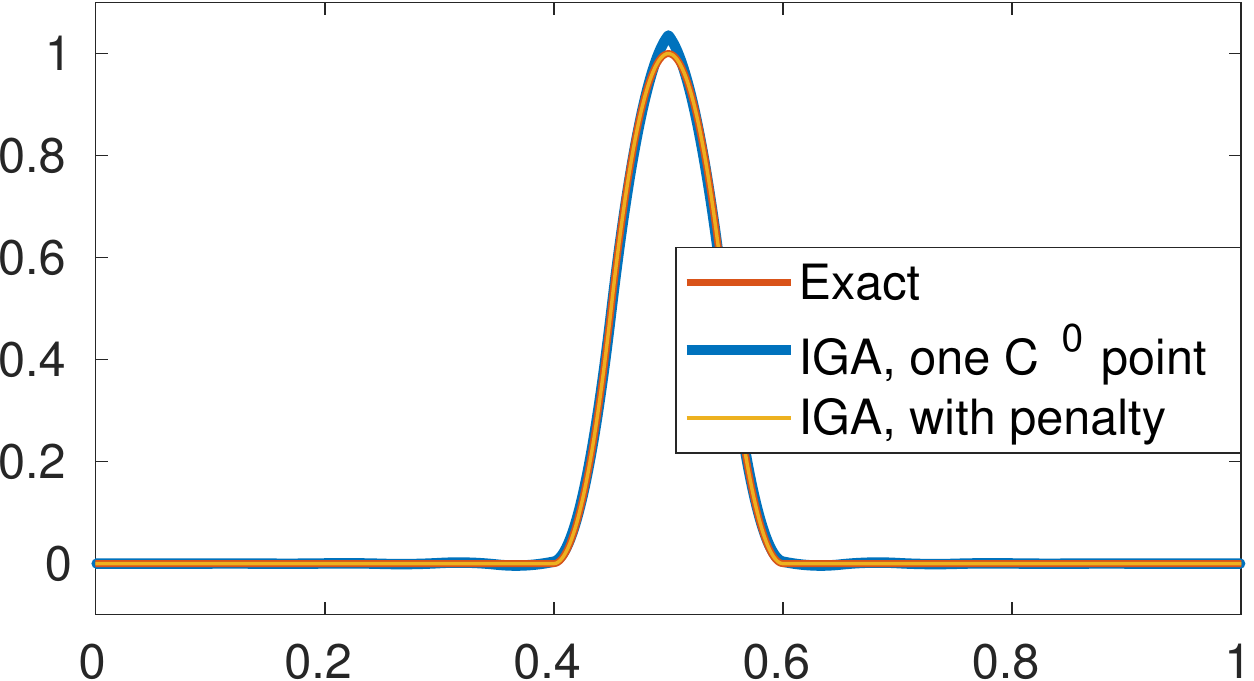}
\includegraphics[width=.49\textwidth]{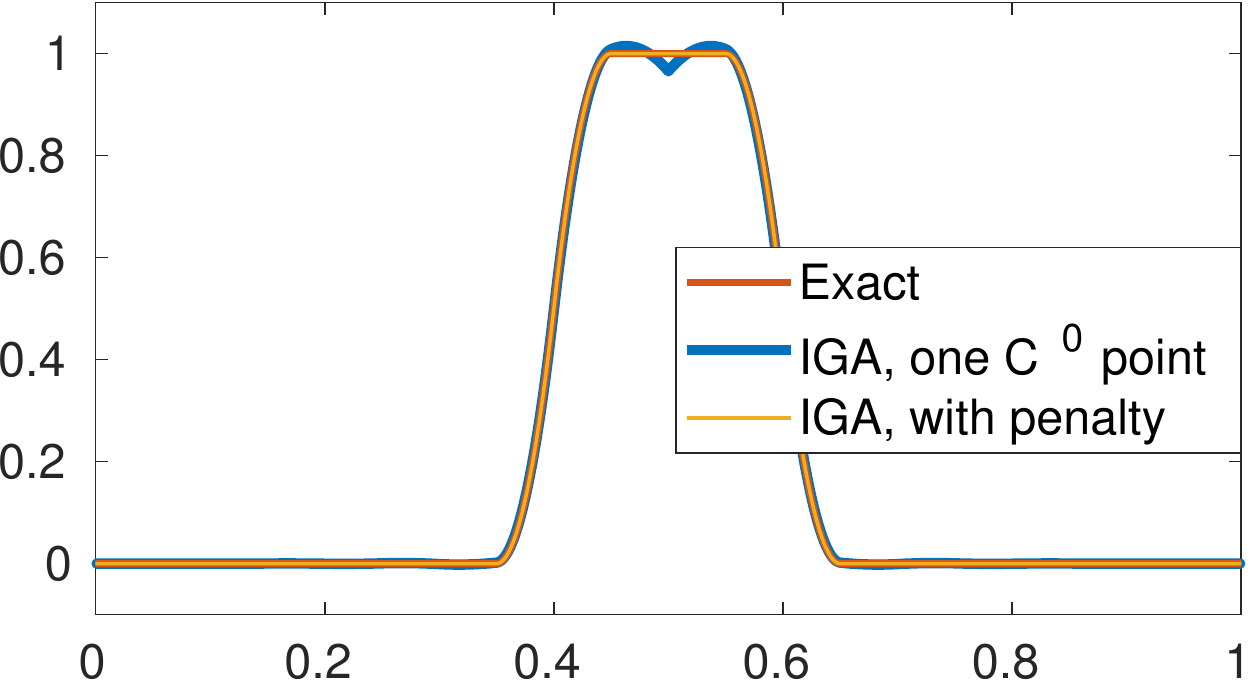}
\caption{$L^2$ best-approximation (penalty value 100).}
\label{fig:ex_simple_approximation}
\end{center}\end{figure}

Let us denote the eigenvectors as $u_{h,i}\in V_h$ for the standard  case with a $C^0$ point and $\widehat u_{h,i} \in V_h$ for the penalty case, with $i=1,\ldots,21$, and note that both  sets span $V_h$. 
Removing the largest eigenvalue by a low-rank modification is equivalent to restricting the solution to the subspace of the first $20$ eigenvalues:
\[
V_h' = \operatorname{span}\{u_{h,i}\in V_h,~ i=1,\ldots,20\}, \quad \widehat V_h' = \operatorname{span}\{\widehat u_{h,i}\in  V_h, ~i=1,\ldots,20\}.
\]
In Figure~\ref{fig:ex_simple_approximation}, we study the best-approximation properties of these reduced spaces by computing the $L^2$-projection of two smooth splines $v_h\in V_h$ onto the spaces $V_h'$ and $\widehat V_h'$. We see a significantly better approximation  in the penalized space $\widehat V_h'$, while we clearly see the non-smoothness of the approximation in $V_h'$. 
This can also be seen in terms of the relative $L^2$ error of the projection, which is presented in Table~\ref{tab:simple_ex_approximation_errors} for different choices of  penalty. With growing values of the penalty parameter, the target splines are approximated more and more precisely.
\begin{table}\begin{center}\footnotesize
\begin{tabular}{r|c|c|c|c|c|c}
penalty value  & 0 &  0.01 & 1 & 10 & 100 &  10\,000 \\ \hline
\vphantom{$\Big($} Example 1 &  0.0185 &  0.0170 &  0.0019  & $2.0590\cdot 10^{-4}$ & $2.0800\cdot 10^{-5}$ &  $2.0823\cdot 10^{-7}$\\
\vphantom{$\Big($} Example 2 & 0.0118 &   0.0107 &   0.0011  & $1.1822\cdot 10^{-4}$ & $1.1928\cdot 10^{-5}$ &  $1.1940\cdot 10^{-7}$
\end{tabular}
\caption{Relative $L^2$ projection error for the two smooth splines shown in Figure~\ref{fig:ex_simple_approximation} for different values of the penalty.}
\label{tab:simple_ex_approximation_errors}
\end{center}\end{table}

This confirms that the imposition of higher regularity through a penalty approach constitutes  a simple way to recover the approximation properties of a $C^1$ space. 
This may have important implications, e.g., in dynamics problems where high modes participate to the response of the analyzed structure. 
It also yields  a better CFL condition, which allows  larger stable time integration steps. In particular in the IGA framework, this is relevant also when a consistent mass is used, since lumped mass is known to be limited to second-order accuracy even for higher orders. This is the reason why there is a strong research interest in predictor-multicorrector explicit algorithms (see, e.g.,~\cite{auricchio:12,evans:18}) making use of the consistent mass for the evaluation of the residual vector and techniques to directly assemble an approximate (banded) inverse of the consistent mass~\cite{hiemstra2017,reali2016}.

\subsection{One-dimensional results} 
In this subsection, we report on one-dimensional results obtained with  higher-order penalty couplings.
On the unit line, the set of eigenvectors with pure Neumann conditions is given by $u_n(x) = \sqrt 2 \sin(n \pi x)$, with the corresponding eigenvalue $\lambda_n = n^2\pi^2$, $n=0,1,\ldots$, see~\cite{hughes:14}.

\begin{figure}[htb]
\begin{center}
\includegraphics[height=8.85em]{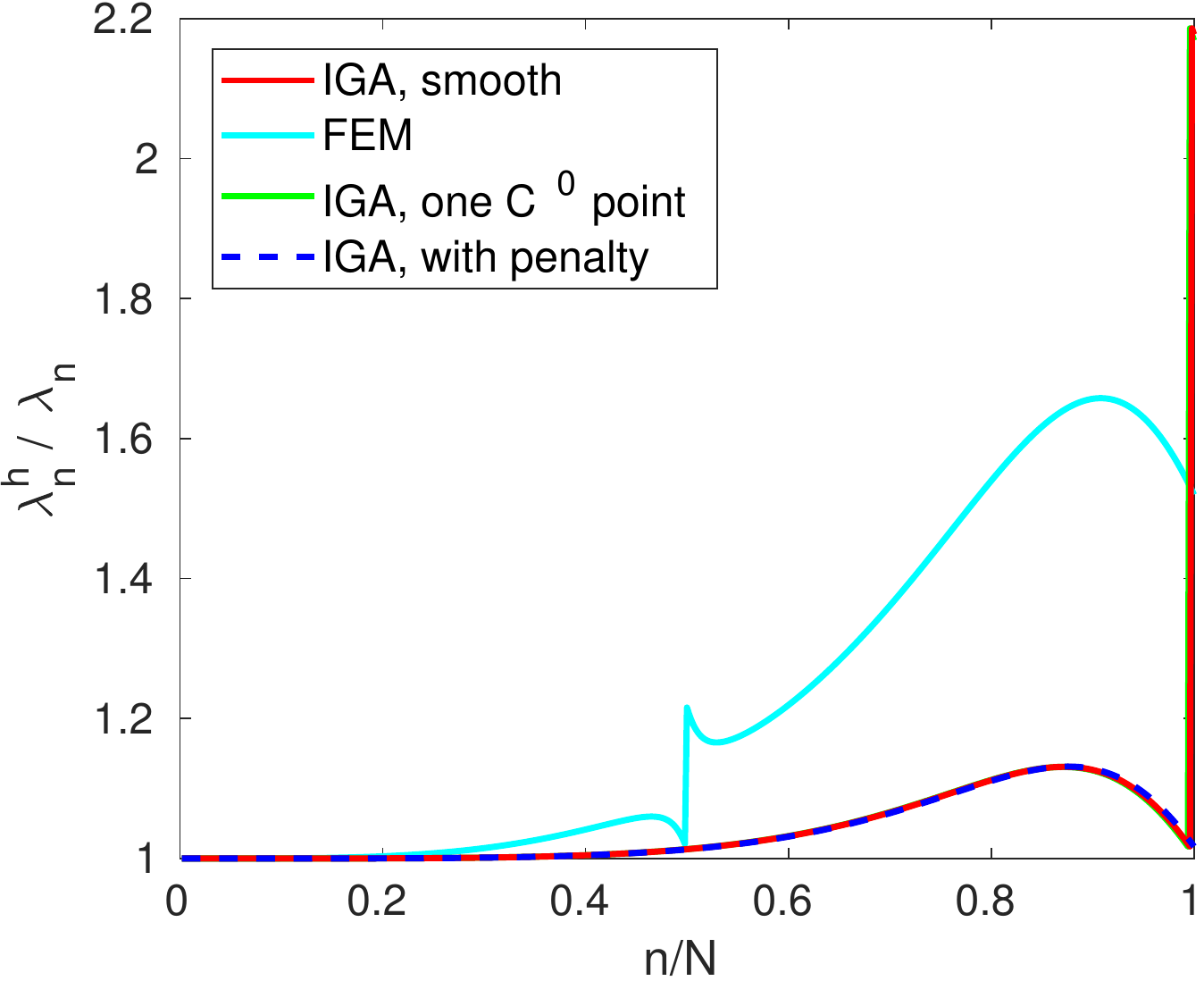}
\includegraphics[height=8.85em]{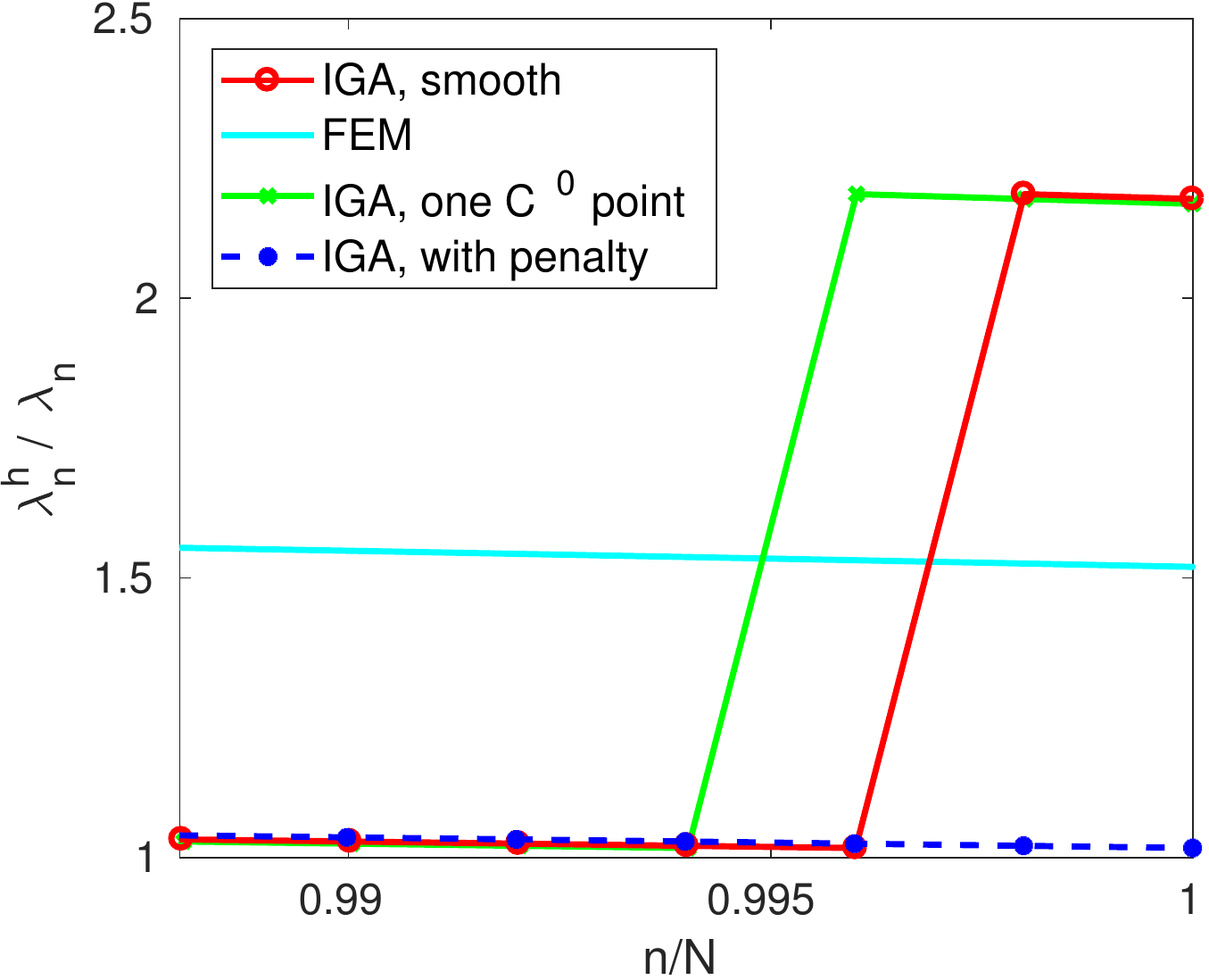}
\includegraphics[height=8.85em]{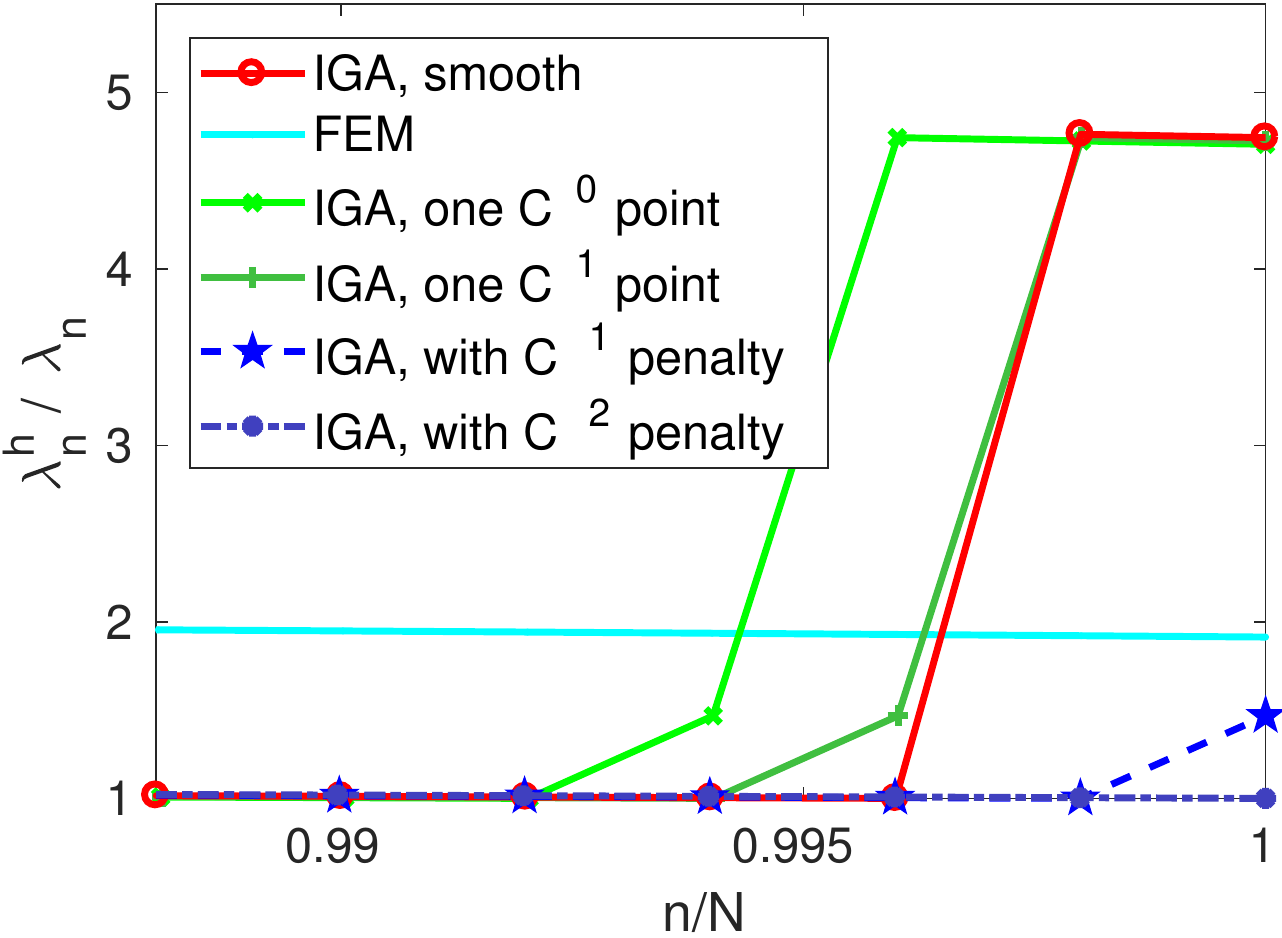}
\caption{One-dimensional discrete spectrum. Left: Entire normalized discrete spectrum for $p=2$. Middle:  Zoom of the last part of the normalized discrete spectrum for $p=2$  Right: Zoom of the last part of the normalized discrete spectrum for $p=3$.}
\label{fig:neum_ord_2_and_3_with_zoom}
\end{center}
\end{figure}

As it is well-known~\cite{hughes:08} finite elements fail to approximate the higher part of the spectrum while IGA with maximal regularity allows to obtain good results. However IGA with reduced regularity, e.g., introduced
by a $C^0$-line or by a weak mortar coupling across multiple-patches, also introduces
outliers at the high frequency end of the spectrum.  By penalizing the jumps in the normal derivatives,
the number of these  outliers can be significantly reduced.
Due to the Neumann boundary conditions, we also see outliers for the smooth space. However, the penalty can reduce even these outliers and we end up with better results than with the original smooth spline space, see Figure~\ref{fig:neum_ord_2_and_3_with_zoom}.  
The results for degrees 4 and 5 are similar to those shown for degree 2 and 3 and, therefore, are not reported here.

\subsection{Application to linear elasticity}
Now, we apply the penalty method to a non-trivial example of elasticity. We reconsider the two-dimensional beam with three circular cut-outs, see the top left of Figure~\ref{fig:beam_with_holes}, clamped on the left side, with Neumann boundaries on the remaining edges including the circular holes.

We solve the eigenvalue problem of linear elasticity:
\[
-\divergence \stress(\mathbf u) = \lambda \mathbf u  \quad \text{ in } \Omega
\]
where the linearized stress and strain are given by
$\stress(\mathbf u) = 2\bar \mu \strain(\mathbf u) + \bar \lambda \operatorname{tr} \strain(\mathbf u)\mathbf I$ and  $\strain(\mathbf u) = (\nabla\mathbf u + \nabla\mathbf u^\top)/2$, respectively. The Lam\'e parameters depend on the elastic modulus $E=1$ and Poisson's ratio $\nu = 0.3$ as 
$\bar \mu = E/( 2+2\nu)$ and $\bar \lambda = \nu E/( (1+\nu)(1-2\nu) )$.
For the equations of elasticity, the surface traction $\stress(\mathbf u) \mathbf n$ plays the role of the normal derivative in the Laplace setting. Hence, the normal derivative in the penalty terms is replaced by $\stress(\mathbf u) \mathbf n$.

We note that for such applications, the different penalty parameters must be well-balanced to ensure a good separation of the physical and unphysical eigenvalues. In this example, it turned out best to only consider the Neumann penalty terms, since the outliers of the Neumann boundary dominate the spectrum. For practical applications, balancing the different penalty terms  can be performed on a coarse mesh with low cost. 

\begin{figure}[ht]
\begin{center}
\begin{minipage}{.47\textwidth}\includegraphics[width=\textwidth]{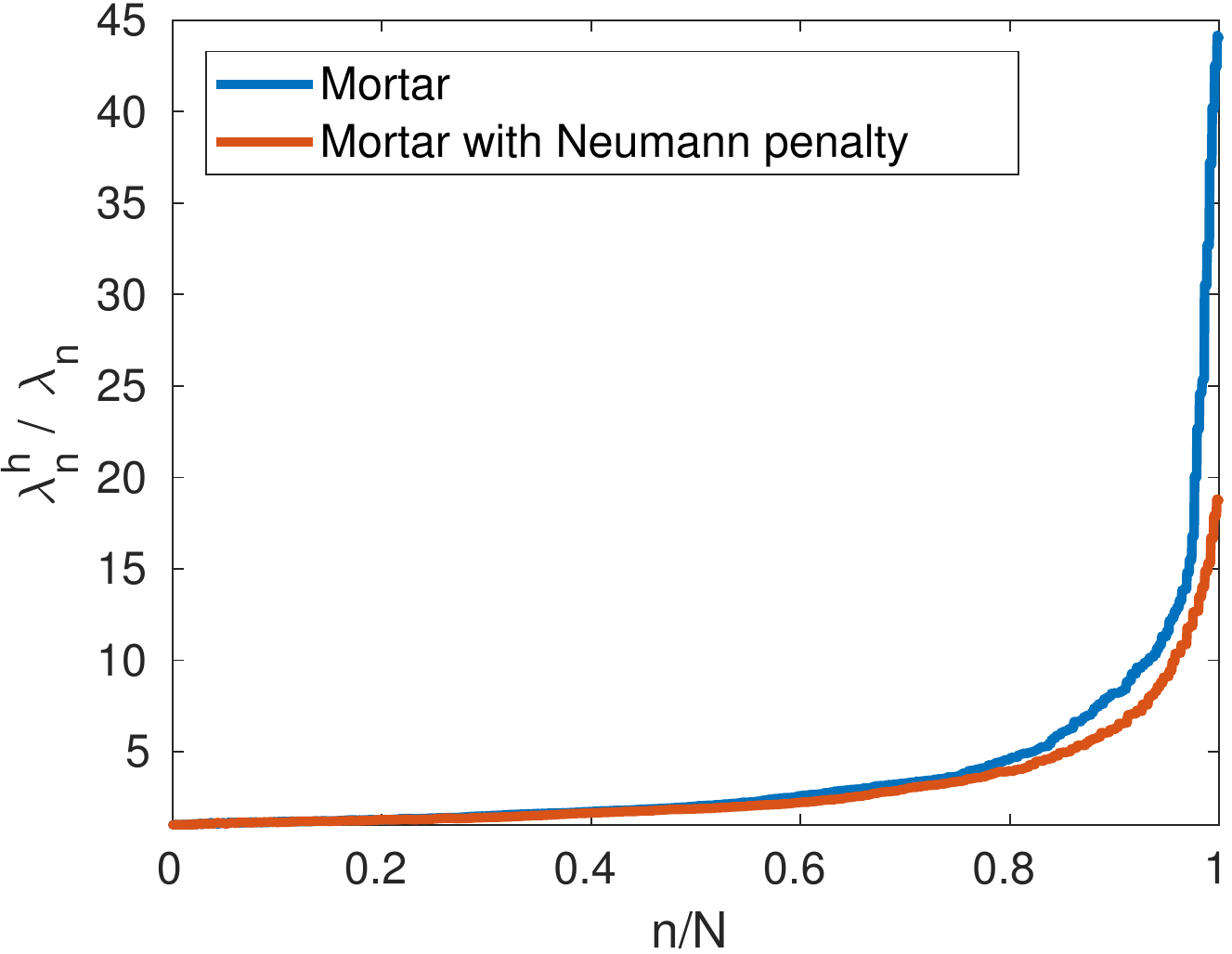}\end{minipage} \hspace{1em}
\begin{minipage}{.47\textwidth}\includegraphics[width=\textwidth]{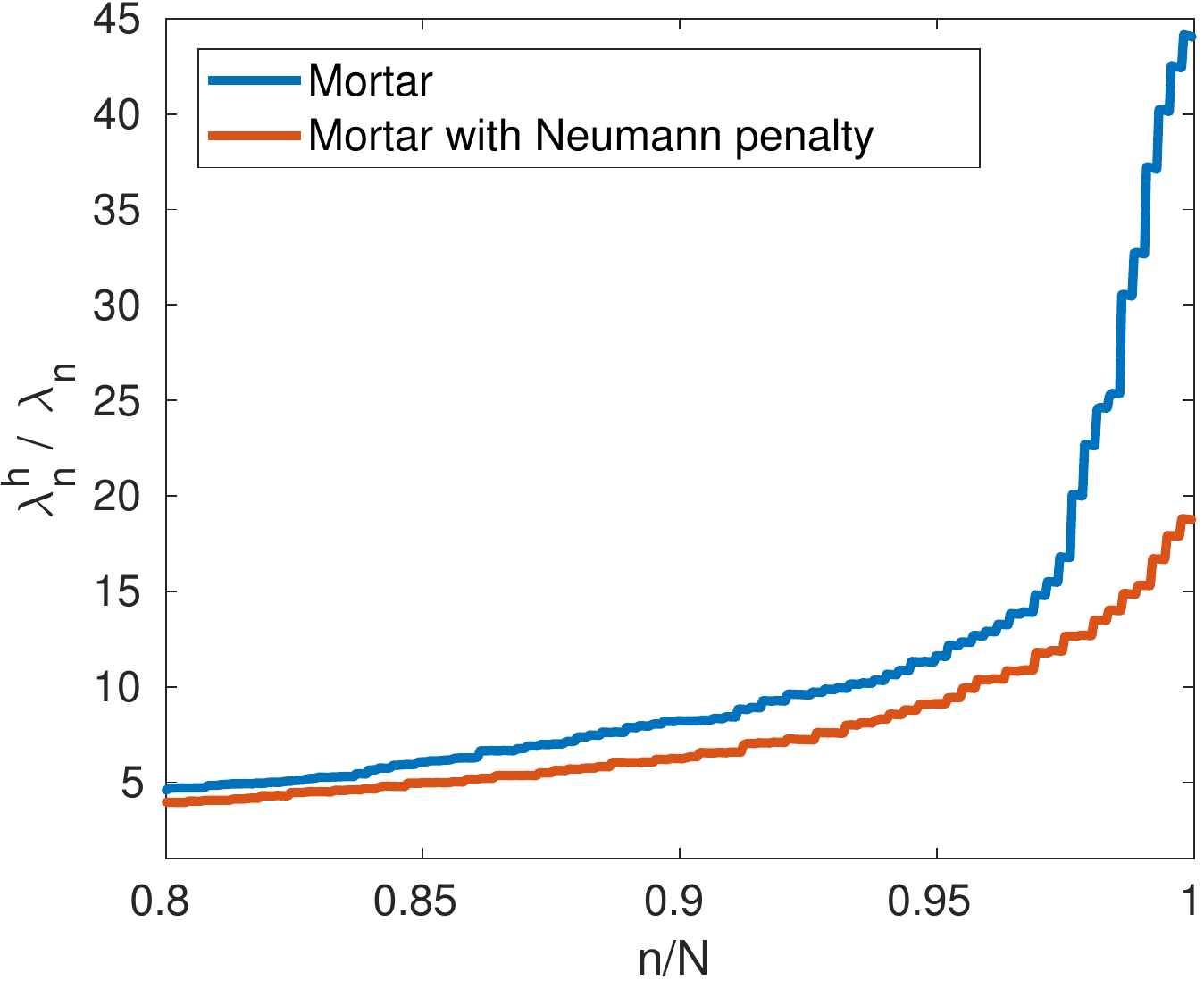}\end{minipage}
\caption{Left: Normalized discrete spectra for the linear elastic beam with circular cut-outs. Right: Zoom to the last 20\% of the spectra. }
\label{fig:laplace_beam}
\end{center}
\end{figure}

As we have no exact solution, we compare the results to a computed reference solution.
The results for a quadratic discretization are shown in Figure~\ref{fig:laplace_beam}. 
Here we have chosen a penalty parameter of $C=10^5$ and note that a large penalty parameter guarantees a clear separation of the physical and unphysical eigenvalues. 
We see that the proposed method provides a significant overall improvement of the discrete spectrum also in the framework of elasticity. In particular, the maximal outlier is reduced to less then half of its value.
Indeed we see that, even though 12 interfaces are present, the Neumann outliers are dominating the spectral approximation and are removed by the proposed penalty.

\section{Vibroacoustical application with a fourth order eigenvalue problem}
\label{sec:numerics_violin}
\begin{figure}[htb!]
\begin{minipage}{.48\textwidth}\begin{center}
\includegraphics[width=.75\textwidth]{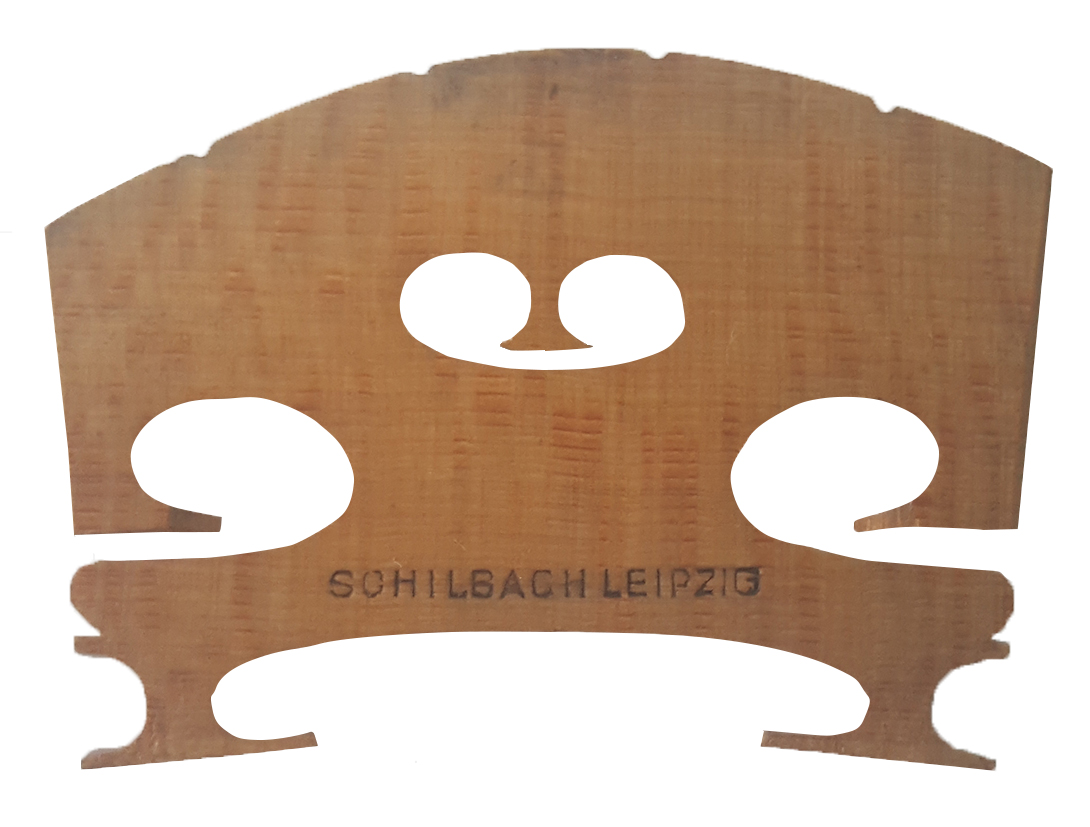}
\caption{Bridge of a violin}
\label{fig:violinsteg_foto}
\end{center}\end{minipage}
\begin{minipage}{.48\textwidth}\begin{center}
\includegraphics[width=.75\textwidth]{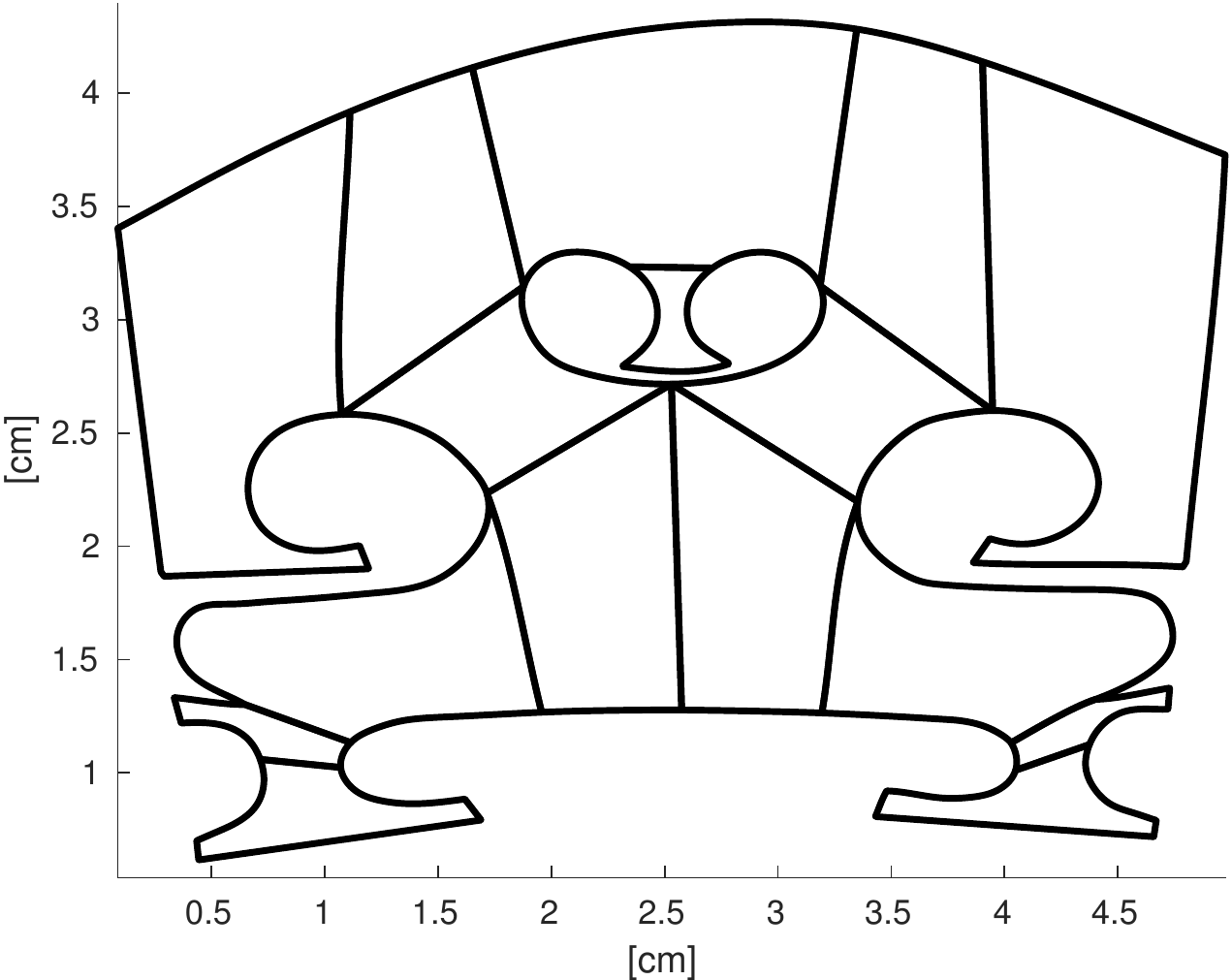}
\caption{Multi-patch representation of the bridge with 16 patches}
\label{fig:violinsteg_decomp}
\end{center}\end{minipage}
\end{figure}
As a final example, let us consider a vibroacoustical example. 
 The bridge of a violin~\cite{horger:17} shown in  Figure~\ref{fig:violinsteg_foto} has an important influence on the acoustics of the instrument. As the geometry is rather thin (thickness of approx. 1\,mm), a plate mode is convenient for an analysis of the out-of-plane eigenmodes. 
In vibroacoustics one is interested in the first part of the spectrum, so we solve the biharmonic eigenvalue problem
\begin{align*}
\Delta\Delta u &= \lambda u\quad \text{ in } \Omega, \\
u=0, \quad &\partial_{\bf n} u = 0  \quad \text { on } \ \GammaD, \\
\partial_{\bf n}^2 u = 0, \quad& \left( \nabla  \Delta u + \boldsymbol \Psi u \right) \cdot \mathbf n = 0 \quad \text{ on } \GammaN,
\end{align*}
with $\boldsymbol \Psi u = ( \partial_x \partial_y^2 u, \partial_y \partial_x^2 u)^\top$ and  homogeneous Dirichlet conditions applied to the bottom of the two `feet' and  natural boundary conditions on the remaining boundary parts.
We consider such isotropic material laws, that the resulting eigenproblem can be rescaled to the one stated above. In this case, changing the elastic modulus   only influences the eigenvalue and not the eigenmode, which allows us to solve the unweighted bilaplace equation. For the more general case, we refer to~\cite{reali:15}. More complex models can also take into account the different behavior of wood in both coordinate directions by considering an  orthotropic Kirchhoff plate.
\begin{figure}
\begin{minipage}{.3\textwidth}\begin{center}
\includegraphics[width=\textwidth]{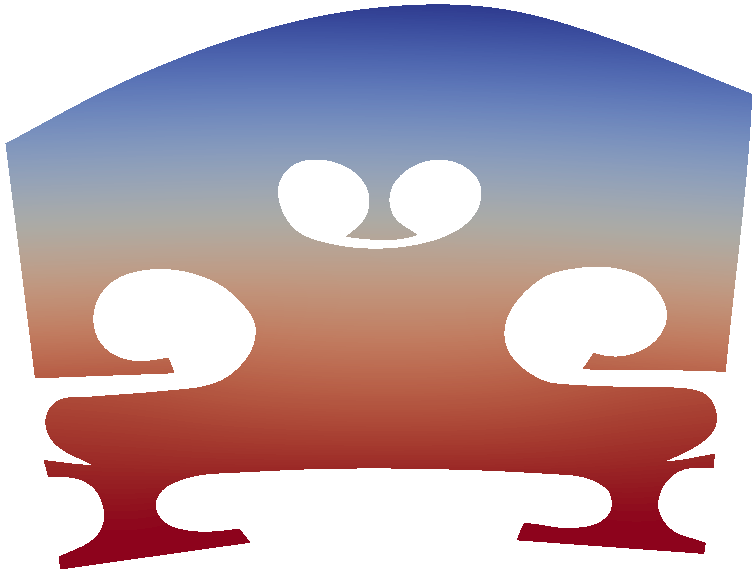}\\
$1^{\rm st}$ eig.value:  0,0221
\end{center}\end{minipage}\hfill
\begin{minipage}{.3\textwidth}\begin{center}
\includegraphics[width=\textwidth]{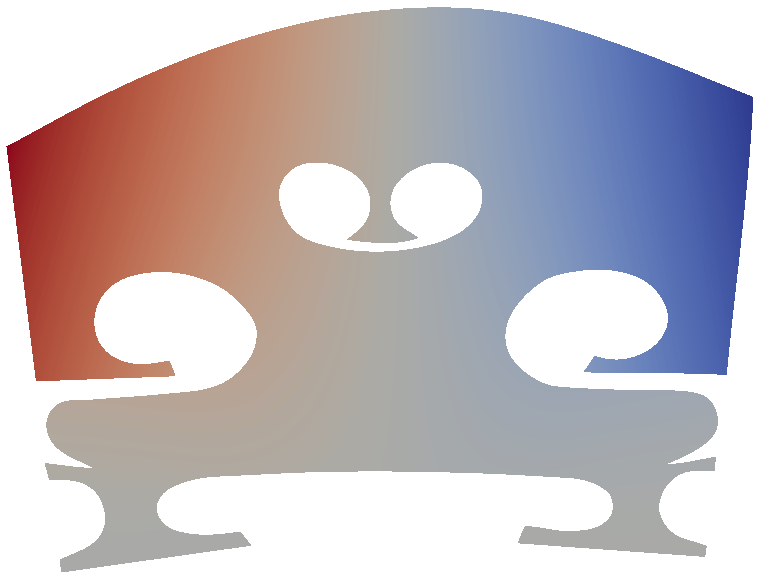}\\
$2^{\rm nd}$ eig.value: 0,136\end{center}\end{minipage}\hfill
\begin{minipage}{.3\textwidth}\begin{center}
\includegraphics[width=\textwidth]{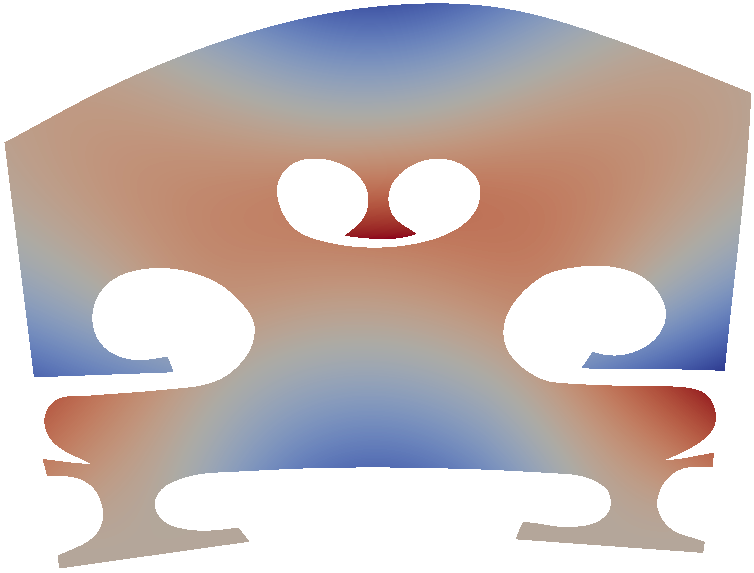}\\
$8^{\rm th}$ eig.value:  9,71
\end{center}\end{minipage}
\caption{Plot of the first, second and eighth eigenmodes with the corresponding eigenvalue}
\label{fig:violin_plot_eigenmodes}
\end{figure}

Since the first part of the spectrum is not influenced by the outliers, we do not use the penalty to improve the high eigenmodes as  for the previous eigenvalue problems. Instead, we use the penalty to be able to solve the fourth order plate problem with the $H^2$-nonconforming mortar space.
Thus, we use the bilinear form $a_h^{\rm{bi}}$ introduced in Section~\ref{sec:consistency_terms}, which includes the first order penalty and consistency terms for the plate problem and solve the following problem. 
Find $ (u_h, \widehat{\tau}_h) \in V_h \times M_h$, $\lambda_h \in \R$, such that 
 \begin{align*}
		a_h^{\rm{bi}}(u_h, v_h)+ b({ \widehat{\tau}_h, v_h }) &= \lambda_h \, m(u_h, v_h), \quad   v_h \in V_h,\\
		b({ \tau_h, u_h }) &= 0, \quad  \tau_h \in M_h,
\end{align*}
The use of a penalty on the normal derivative to solve the plate eigenvalue problem is also applied in a FEM context by the $C^0$-IPDG method~\cite{brenner:15}.

The geometry is represented by 16 patches coupled across 16 interfaces as shown in Figure~\ref{fig:violinsteg_decomp}. However since some patches have corners, there are $C^0$-lines within some patches, where the penalty coupling needs to be applied as well, yielding a total of 28 interfaces for the penalty coupling, where we chose  $C=10$. 

\begin{figure}[htb!]
 \includegraphics[width=.49\textwidth]{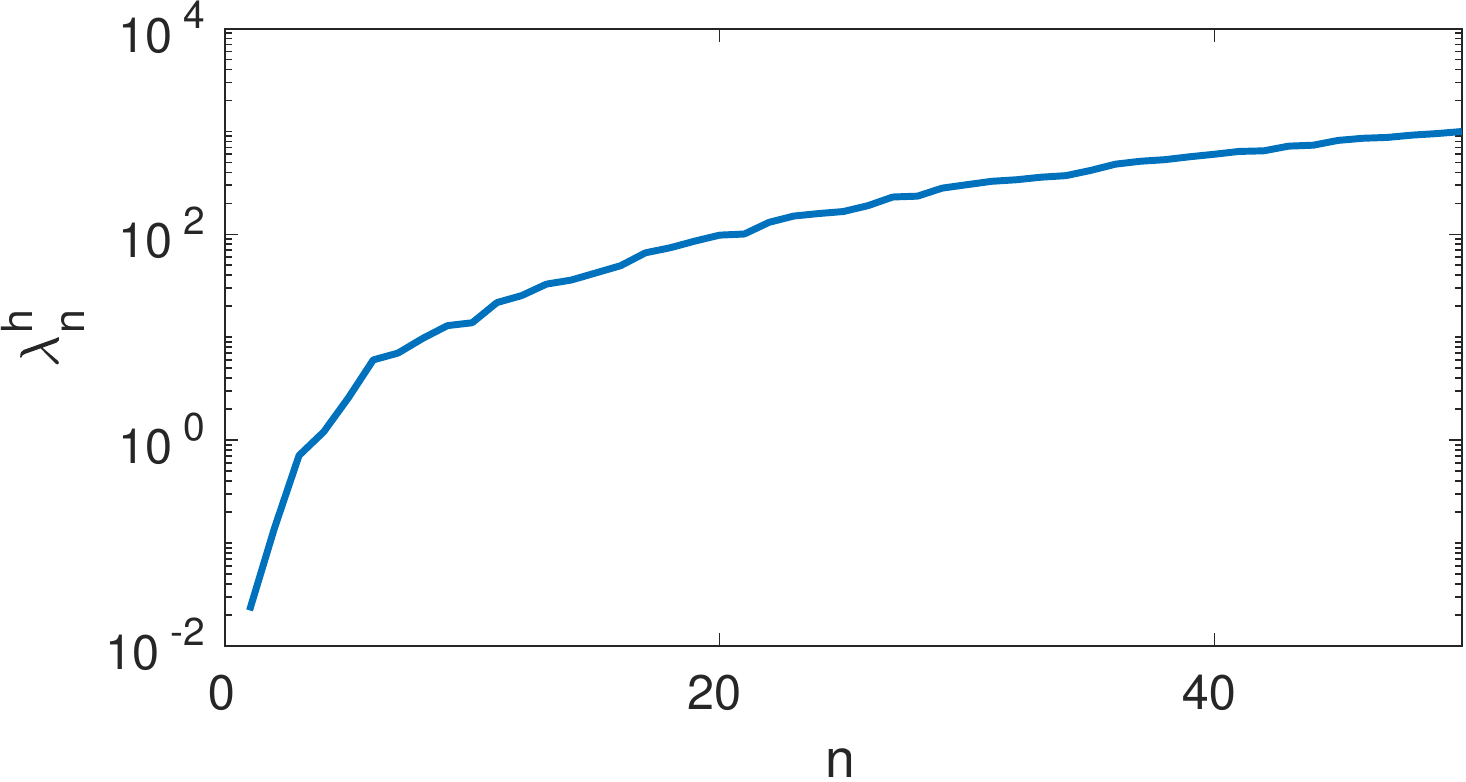} 
\includegraphics[width=.49\textwidth]{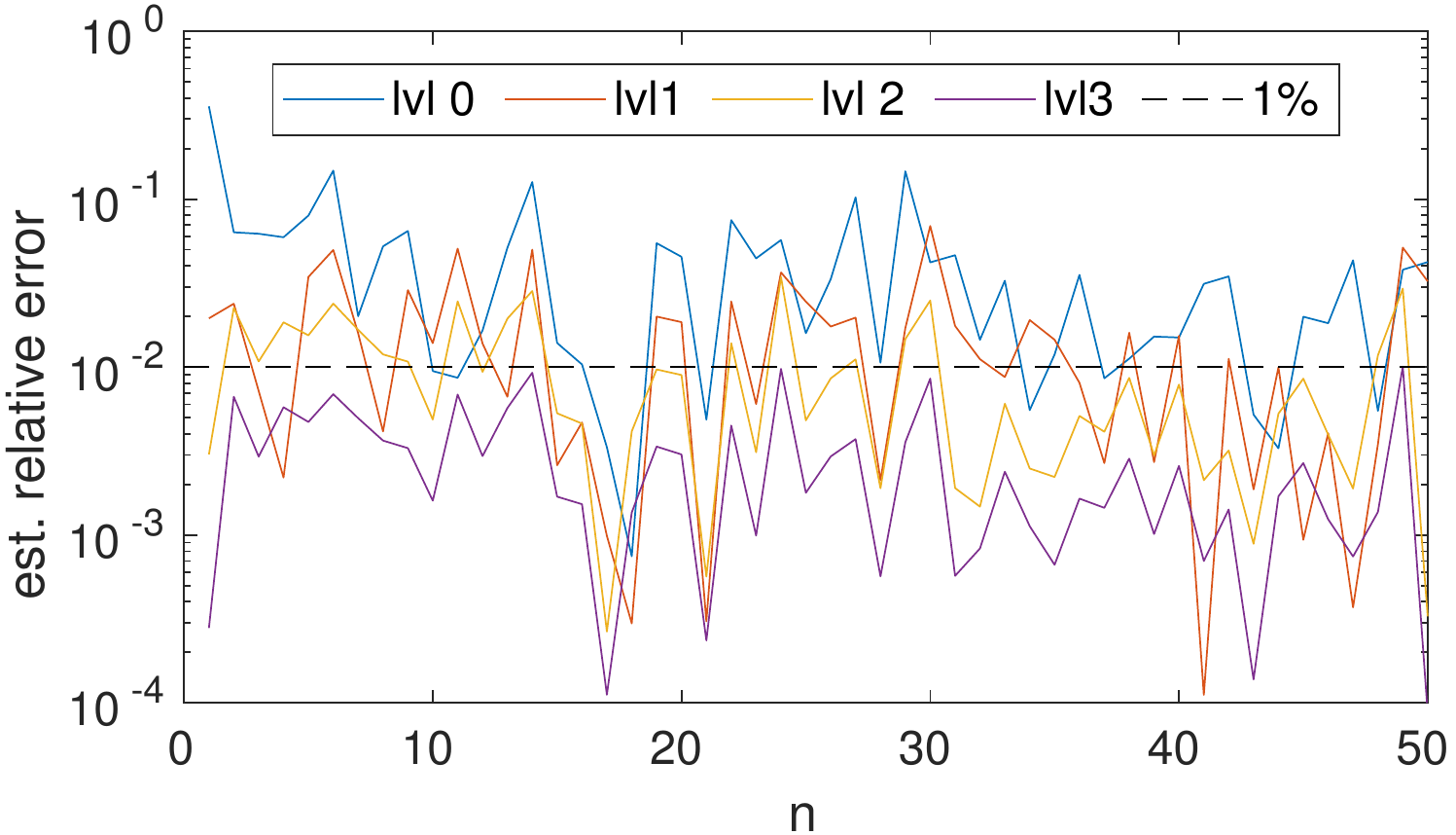}
\caption{First 50 eigenvalues of the biharmonic equation for the violin bridge. Left: discrete eigenvalues on the reference mesh level 4. Right: estimated error values.}
\label{fig:violin_plate_50_eigenvalues}
\end{figure}

A selection of eigenmodes and the corresponding eigenvalues on mesh level~3 for $p=3$ with $33{,}440$ degrees of freedom are shown in Figure~\ref{fig:violin_plot_eigenmodes}. In all cases we see smooth results thanks to the hybrid coupling and in particular no spurious oscillations are observed. 
The first 50 eigenvalues on several mesh levels as well as an estimated error are shown in Figure~\ref{fig:violin_plate_50_eigenvalues}. Here, we see a very good approximation of the relevant eigenmodes for vibroacoustics already on the first meshes. On the finest mesh, level 3, the relative error of all first 50 eigenvalues is below 1\%.

\section{Conclusions}
\label{sec:conclusion}
In this paper, we have studied, in the framework of isogeometric analysis, the effects of higher-order penalty terms for multi-patch geometries and Neumann boundaries on second and fourth order partial differential equations. 
In the context of fourth order problems, the hybrid coupling poses a flexible discretization for multi-patch geometries and can include complicated boundary conditions. 
For second order eigenvalue problems, the hybrid coupling reduces so-called outlier eigenvalues, which is relevant  in several applications such as, e.g., explicit dynamics.

\bibliography{lit}

\begin{thebibliography}{10}
\expandafter\ifx\csname url\endcsname\relax
  \def\url#1{\texttt{#1}}\fi
\expandafter\ifx\csname urlprefix\endcsname\relax\def\urlprefix{URL }\fi
\expandafter\ifx\csname href\endcsname\relax
  \def\href#1#2{#2} \def\path#1{#1}\fi

\bibitem{hughes:05}
T.~J.~R. Hughes, J.~A. Cottrell, Y.~Bazilevs, {Isogeometric analysis: CAD,
  finite elements, NURBS, exact geometry and mesh refinement}, Comput. Methods.
  Appl. Mech. Eng. 194 (2005) 4135--4195.

\bibitem{Hoellig03}
K.~H\"ollig, {Finite Element Methods with B-Splines}, Frontiers in Applied
  Mathematics, SIAM, 2003.

\bibitem{beirao:14}
L.~Beir{\~a}o~da Veiga, A.~Buffa, G.~Sangalli, R.~V\'{a}zquez, Mathematical
  analysis of variational isogeometric methods, Acta Numer. 23 ({2014})
  157--287.

\bibitem{nguyen:15}
V.~P. Nguyen, C.~Anitescu, S.~P. Bordas, T.~Rabczuk, Isogeometric analysis: An
  overview and computer implementation aspects, Math. Comp. Simul. 117 (2015)
  89 -- 116.

\bibitem{ben_belgacem:99}
F.~{Ben Belgacem}, {The mortar finite element method with Lagrange
  multipliers}, Numer. Math. 84 (1999) 173--197.

\bibitem{bernardi:94}
C.~Bernardi, Y.~Maday, A.~T. Patera, A new nonconforming approach to domain
  decomposition: the mortar element method, in: H.~B. et.al. (Ed.), Nonlinear
  partial differential equations and their applications., Vol.~XI, Coll\`{e}ge
  de France, 1994, pp. 13--51.

\bibitem{Woh01_a}
B.~Wohlmuth, Discretization Techniques and Iterative Solvers Based on Domain
  Decomposition, Vol.~17, Springer, Heidelberg, 2001.

\bibitem{hesch:12}
C.~Hesch, P.~Betsch, Isogeometric analysis and domain decomposition methods,
  {Comput. Methods Appl. Mech. Eng.} 213--216 (2012) 104--112.

\bibitem{bletzinger:14}
A.~Apostolatos, R.~Schmidt, R.~W{\"u}chner, K.-U. Bletzinger, {A Nitsche-type
  formulation and comparison of the most common domain decomposition methods in
  isogeometric analysis}, Int. J. Numer. Methods Eng. 97 (2014) 473--504.

\bibitem{dornisch:14}
W.~Dornisch, G.~Vitucci, S.~Klinkel, The weak substitution method – an
  application of the mortar method for patch coupling in {NURBS}-based
  isogeometric analysis, Int. J. Numer. Methods in Eng. 103~(3) (2015)
  205--234.

\bibitem{nguyen:14}
V.~P. Nguyen, P.~Kerfriden, M.~Brino, S.~P.~A. Bordas, E.~Bonisoli, Nitsche's
  method for two and three dimensional {NURBS} patch coupling, Comput. Mech.
  53~(6) (2014) 1163--1182.

\bibitem{popp:18}
L.~Wunderlich, A.~Seitz, M.~D. Alaydin, B.~Wohlmuth, A.~Popp, Biorthogonal
  splines for optimal weak patch-coupling in isogeometric analysis with
  applications to finite deformation elasticity,
  \texttt{https://arxiv.org/abs/1806.11535}.

\bibitem{brivadis:15}
E.~Brivadis, A.~Buffa, B.~Wohlmuth, L.~Wunderlich, Isogeometric mortar methods,
  Comput. Methods Appl. Mech. Eng. 284 (2015) 292--319.

\bibitem{kiendl:09}
J.~Kiendl, K.-U. Bletzinger, J.~Linhard, R.~W{\"u}chner, Isogeometric shell
  analysis with {Kirchhoff-Love} elements, Comput. Methods Appl. Mech. Eng. 198
  (2009) 3902 -- 3914.

\bibitem{kiendl:10}
J.~Kiendl, Y.~Bazilevs, M.-C. Hsu, R.~W{\"u}chner, K.-U. Bletzinger, The
  bending strip method for isogeometric analysis of {Kirchhoff-Love} shell
  structures comprised of multiple patches, Comput. Methods Appl. Mech. Eng.
  199 (2010) 2403 -- 2416.

\bibitem{bouclier:16}
R.~Bouclier, J.-C. Passieux, M.~Sala\"un, Development of a new, more regular,
  mortar method for the coupling of {NURBS} subdomains within a {NURBS} patch:
  Application to a non-intrusive local enrichment of {NURBS} patches, Comput.
  Methods Appl. Mech. Eng. 316 (2017) 123--150.

\bibitem{collin:16}
A.~Collin, G.~Sangalli, T.~Takacs, Analysis-suitable {$G^1$} multi-patch
  parametrizations for {$C^1$} isogeometric spaces, Comput. Aided Geom. Design
  47 (2016) 93 -- 113.

\bibitem{coox:16}
L.~Coox, F.~Greco, O.~Atak, D.~Vandepitte, W.~Desmet, A robust patch coupling
  method for {NURBS}-based isogeometric analysis of non-conforming multipatch
  surfaces, Comput. Methods Appl. Mech. Eng. 316 (2017) 235--260.

\bibitem{baker:77}
G.~A. Baker, Finite element methods for elliptic equations using nonconforming
  elements, Math. Comp. 31~(137) (1977) 45--59.

\bibitem{mozolevski:03}
I.~Mozolevski, E.~S{\"u}li, A priori error analysis for the hp-version of the
  discontinuous {Galerkin} finite element method for the biharmonic equation,
  Comput. Methods Appl. Mech. Eng. 3 (2003) 596--607.

\bibitem{cottrell:06}
J.~A. Cottrell, A.~Reali, Y.~Bazilevs, T.~J.~R. Hughes, Isogeometric analysis
  of structural vibrations, Comput. Methods Appl. Mech. Eng. 195~(41--43)
  (2006) 5257 -- 5296.

\bibitem{hughes:08}
T.~J.~R. Hughes, A.~Reali, G.~Sangalli, Duality and unified analysis of
  discrete approximations in structural dynamics and wave propagation:
  Comparison of p-method finite elements with k-method {NURBS}, Comput. Methods
  Appl. Mech. Eng. 197~(49–50) (2008) 4104 -- 4124.

\bibitem{hughes:09}
J.~A. Cottrell, T.~J.~R. Hughes, Y.~Bazilevs, Isogeometric Analysis. Towards
  Integration of CAD and FEA, Wiley, Chichester, 2009.

\bibitem{hughes:14}
T.~J.~R. Hughes, J.~A. Evans, A.~Reali, Finite element and {NURBS}
  approximations of eigenvalue, boundary-value, and initial-value problems,
  Comput. Methods Appl. Mech. Eng. 272 (2014) 290 -- 320.

\bibitem{calo:18}
V.~Puzyrev, Q.~Deng, V.~Calo, Spectral approximation properties of isogeometric
  analysis with variable continuity, Comput. Methods Appl. Mech. Eng. 334
  (2018) 22 -- 39.

\bibitem{takacs:16}
S.~Takacs, T.~Takacs, Approximation error estimates and inverse inequalities
  for {B}-splines of maximum smoothness, Math. Models Methods Appl. Sci.
  26~(07) (2016) 1411--1445.

\bibitem{gallistl:17}
D.~Gallistl, P.~Huber, D.~Peterseim, On the stability of the {Rayleigh}--{Ritz}
  method for eigenvalues, Numer. Math. 137~(2) (2017) 339--351.

\bibitem{hiemstra2017}
R.~R. Hiemstra, A.~Reali, G.~Sangalli, M.~Tani, J.~A. Evans, T.~J.~R. Hughes,
  Efficient isogeometric collocation for explicit structural dynamics:
  High-order mass lumping and outlier removal, in preparation.

\bibitem{reali2016}
A.~Reali, T.~J.~R. Hughes, {IGA} collocation, aka ``the ultimate reduced
  quadrature {IGA} method": Some results, applications, and open problems'',
  in: WCCM XII \& APCOM VI - 12th World Congress on Computational Mechanics and
  6th Asian Pacific Congress on Computational Mechanics, Seoul, 2016.

\bibitem{bazilevs:06}
Y.~Bazilevs, L.~Beir{\~a}o~da Veiga, J.~A. Cottrell, T.~J.~R. Hughes,
  G.~Sangalli, Isogeometric analysis: Approximation, stability and error
  estimates for h-refined meshes, Math. Models Methods Appl. Sci. 16~(7) (2006)
  1031--1090.

\bibitem{piegl:97}
L.~Piegl, W.~Tiller, The NURBS Book, Springer, 1997.

\bibitem{Schumaker:07}
L.~Schumaker, Spline Functions: Basic Theory, 3rd Edition, Cambridge University
  Press, Cambridge, 2007.

\bibitem{grisvard:11}
P.~Grisvard, Elliptic Problems in Nonsmooth Domains, SIAM, Philadelphia, 2011.

\bibitem{brenner:05}
S.~C. Brenner, L.-Y. Sung, {$C^0$} interior penalty methods for fourth order
  elliptic boundary value problems on polygonal domains, J. Sci. Comput. 22~(1)
  (2005) 83--118.

\bibitem{brenner:15}
S.~C. Brenner, P.~Monk, J.~Sun, {$C^0$} interior penalty {Galerkin} method for
  biharmonic eigenvalue problems, in: R.~M. Kirby, M.~Berzins, J.~S. Hesthaven
  (Eds.), Spectral and High Order Methods for Partial Differential Equations
  ICOSAHOM 2014, Springer International Publishing, Cham, 2015, pp. 3--15.

\bibitem{geopdes:11}
C.~de~Falco, A.~Reali, R.~V\'{a}zquez, {GeoPDEs}: A research tool for
  isogeometric analysis of {PDEs}, Adv. Eng. Softw. 42(12) (2011) 1020--1034.

\bibitem{geopdes:16}
R.~V{\'a}zquez, {A new design for the implementation of isogeometric analysis
  in Octave and Matlab: GeoPDEs 3.0}, Comp. Math. Appl. 72~(3) (2016) 523 --
  554.

\bibitem{auricchio:12}
F.~Auricchio, L.~Beir{\~a}o~da Veiga, T.~J.~R. Hughes, A.~Reali, G.~Sangalli,
  Isogeometric collocation for elastostatics and explicit dynamics, Comput.
  Methods Appl. Mech. Eng. 249-252 (2012) 2 -- 14.

\bibitem{evans:18}
J.~A. Evans, R.~R. Hiemstra, T.~J.~R. Hughes, A.~Reali, Explicit higher-order
  accurate isogeometric collocation methods for structural dynamics, Comput.
  Methods Appl. Mech. Eng. 338 (2018) 208 -- 240.

\bibitem{horger:17}
T.~Horger, B.~Wohlmuth, L.~Wunderlich, Reduced basis isogeometric mortar
  approximations for eigenvalue problems in vibroacoustics, in: P.~Benner,
  M.~Ohlberger, A.~Patera, G.~Rozza, K.~Urban (Eds.), Model Reduction of
  Parametrized Systems, Springer International Publishing, Cham, 2017, pp.
  91--106.

\bibitem{reali:15}
A.~Reali, H.~Gomez, An isogeometric collocation approach for
  {Bernoulli}--{Euler} beams and {Kirchhoff} plates, Comput. Methods Appl.
  Mech. Eng. 284 (2015) 623 -- 636.

\end{thebibliography}
\end{document}